\documentclass{amsart} 
\usepackage{amsxtra}
\usepackage{amssymb}
\usepackage{tikz-cd}
\usetikzlibrary{cd,shapes}
\usetikzlibrary{arrows}
\usepackage[mathscr]{eucal}

\input xy \xyoption{matrix} \xyoption{arrow}\xyoption{frame}

\theoremstyle{plain}
\tikzstyle{mutable}=[align=center, draw,inner sep=0pt,shape=ellipse,minimum height=0.4cm, minimum width=0.4cm]
\tikzstyle{frozen}=[inner sep=0.5mm,rectangle,draw]

\newcommand{\nn}{\hfill\nonumber}
\newtheorem{theorem}{Theorem}[section]
\newtheorem{lemma}[theorem]{Lemma}
\newtheorem{definition-theorem}[theorem]{Definition-Theorem}
\newtheorem{proposition}[theorem]{Proposition}
\newtheorem{corollary}[theorem]{Corollary}
\newtheorem{definition}[theorem]{Definition}
\newtheorem{remark}[theorem]{Remark}
\newtheorem{conjecture}[theorem]{Conjecture}
\newtheorem{notation}[theorem]{Notation}

\newtheorem*{maintheorem*}{Main Theorem}
\newtheorem*{theorem*}{Theorem}
\newtheorem*{theoremA*}{Theorem A}
\newtheorem*{theoremB*}{Theorem B}
\theoremstyle{definition}

\newtheorem*{definition*}{Definition}
\newcommand \bth[1] { \begin{theorem}\label{t#1} }
\newcommand \ble[1] { \begin{lemma}\label{l#1} }

\newcommand \bpr[1] { \begin{proposition}\label{p#1} }
\newcommand \bco[1] { \begin{corollary}\label{c#1} }
\newcommand \bde[1] { \begin{definition}\label{d#1}\rm }
\newcommand \bex[1] { \begin{example}\label{e#1}\rm }
\newcommand \bre[1] { \begin{remark}\label{r#1}\rm }
\newcommand \bcj[1] { \begin{conjecture}\label{j#1}\rm }

\newcommand \bnota[1] { \begin{notation}\label{n#1}\rm }
\renewcommand {\eth} { \end{theorem} }
\newcommand {\ele} { \end{lemma} }

\newcommand {\epr} { \end{proposition} }
\newcommand {\eco} { \end{corollary} }
\newcommand {\ede} { \end{definition} }
\newcommand {\eex} { \end{example} }
\newcommand {\ere} { \end{remark} }
\newcommand {\ecj} { \end{conjecture} }

\newcommand {\enota} { \end{notation} }

\def \KK {{\mathbb K}}

\def \Zset {{\mathbb Z}}

\def \QQ {{\mathcal{Q}}}

\def \PP {{\mathcal{P}}}

\def \RR {{\mathcal{R}}}

\def \LL {{\mathcal{L}}}

\def \al {\alpha}
\def \be {\beta}
\def \vpi {\varpi}

\def \om {\omega}

\def \ga {\gamma}

\def \vp {\varphi}
\def \ep {\epsilon}
\def \De {\Delta}

\def \mt  {\mapsto}

\def \hra {\hookrightarrow}

\def \ci  {\circ}

\def \rcor {\rangle}
\def \lcor {\langle}

\def \id { {\mathrm{id}} }

\def \g  {\mathfrak{g}}   

\DeclareMathOperator \Span { {\mathrm{Span}} }

\DeclareMathOperator \sign {{\mathrm{sign}}}

\DeclareMathOperator \opp { {\mathrm{op}} }

\DeclareMathOperator{\up}{up}

\usepackage[normalem]{ulem}
\usepackage{amsmath,amsthm,amsfonts,amssymb, tikz-cd} 
\usepackage{stackrel}

\usepackage{verbatim}      	
\usepackage{listings} 
\usepackage{mathtools}
\usepackage{epsfig}         	
\usepackage{url}		
\usepackage{epic,eepic,latexsym}
\usepackage{graphicx}
\usepackage{xcolor}
\usepackage{amsxtra}

\usepackage{lcd}
\usepackage{mathrsfs}
\usepackage[centering,text={15.5cm,22cm},
        marginparwidth=20mm,dvips]{geometry}
\usepackage[linktocpage]{hyperref}
\usepackage{lscape}
\usepackage{tabularx}
\usepackage{marginnote}
\usepackage[all]{xy}
\usepackage{enumitem}
\usepackage{stmaryrd}
\usepackage{tipa}
\usepackage{upgreek}

\usetikzlibrary {arrows.meta,positioning}
\usepackage{caption}
\usepackage{subcaption}

\DeclareMathOperator{\tw}{tw}

\DeclareMathOperator{\can}{can}
\DeclareMathOperator{\N}{N}

\def\!{\mskip-\thinmuskip}
\let\?=\overline
\let\:=\colon

\newcommand {\kk} {\Bbbk}

\newcommand {\mm} {{\bf m}}

\newcommand {\sd} {{\bf s}}

\theoremstyle{plain}
 \ifdefined\thm
 \else
   \theoremstyle{plain}
    \newtheorem{thm}[theorem]{Theorem}
 \fi

 \ifdefined\prop
 \else
   \theoremstyle{plain}
    \newtheorem{prop}[theorem]{Proposition}
 \fi\ifdefined\cor
 \else
   \theoremstyle{plain}
    \newtheorem{cor}[theorem]{Corollary}
 \fi

 \ifdefined\lem
 \else
   \theoremstyle{plain}
    \newtheorem{lem}[theorem]{Lemma}
 \fi

 \ifdefined\def
 \else
   \theoremstyle{definition}
   \newtheorem{def}[theorem]{Definition}
 \fi
 
 \ifdefined\defn
 \else
   \theoremstyle{definition}
   \newtheorem{defn}[theorem]{Definition}
 \fi
 
 \ifdefined\example
 \else
   \theoremstyle{definition}
     \newtheorem{example}[theorem]{Example}
 \fi
 
 \ifdefined\eg
 \else
     \theoremstyle{definition}
         \newtheorem{eg}[theorem]{Example}
 \fi
 
 \ifdefined\rem
 \else
     \theoremstyle{remark}
         \newtheorem{rem}[theorem]{Remark}
 \fi
 
 \ifdefined\quest
 \else
     \theoremstyle{plain}
     \newtheorem{quest}[theorem]{Question}
 \fi
 
 \ifdefined\assumption*
 \else
     \theoremstyle{plain}
     \newtheorem*{assumption*}{Assumption}
 \fi

  \newtheorem{conj}[theorem]{Conjecture}

\newcommand{\Q}{{\mathbb Q}}
\newcommand{\R}{{\mathbb R}}
\newcommand{\Z}{{\mathbb Z}}

\newcommand{\frakg}{\mathfrak{g}}

\newcommand{\supp}{\operatorname{supp}}

\usepackage{tikz-cd}
\usepackage{bm}

\newcounter{listequation}

\numberwithin{equation}{section}

\newcommand{\BZ}{\mathrm{BZ}}
\newcommand{\GY}{\mathrm{GY}}

\newcommand{\factor}{\mathsf{p}}

\newcommand{\frGp}{\mathcal{FM}}

\newcommand{\uc}{\underline{c}}

\renewcommand{\cor}{\mathrm{cor}}
\newcommand{\wt}{\mathrm{wt}}

\newcommand{\rev}{\mathrm{rev}}
\newcommand{\lex}{\mathrm{lex}}

\newcommand{\dI}{\dot{I}}
\newcommand{\ddI}{\ddot{I}}
\newcommand{\ddB}{\ddot{B}}
\newcommand{\dB}{\dot{B}}

\newcommand{\cF}{\mathcal{F}}
\newcommand{\Br}{\mathsf{Br}}

\newcommand{\tB}{\widetilde{B}}

\newcommand{\cZ}{\mathcal{Z}}
\newcommand{\cS}{\mathcal{S}}
\newcommand{\LP}{\mathcal{LP}}
\newcommand{\bLP}{\overline{\LP}}

\renewcommand{\N}{\mathbb{N}}

\newcommand{\std}{\mathbf{M}}
\newcommand{\bStd}{\overline{\std}}

\renewcommand{\can}{\mathbf{L}}

\newcommand{\rTri}{\dot{\can}}
\newcommand{\dTri}{\ddot{\can}}
\newcommand{\tTri}{\widetilde{\can}}

\newcommand{\PBW}{\mathbf{X}}
\newcommand{\dPBW}{\mathbf{X}^*}

\newcommand{\dCan}{\mathbf{B}^*}

\newcommand{\K}{\mathbb{K}}
\newcommand{\dkappa}{\ddot{\kappa}}
\newcommand{\rkappa}{\dot{\kappa}}

\newcommand{\dsd}{\ddot{\sd}} 
\newcommand{\rsd}{\dot{\sd}}
\newcommand{\ddiag}{\symm'}
\newcommand{\symm}{\mathsf{d}}
\newcommand{\bi}{{\bf i}}
\newcommand{\ubi}{\underline{\bf i}}

\newcommand{\ubj}{\underline{\bf j}}
\newcommand{\ueta}{\underline{\eta}}
\newcommand{\uzeta}{\underline{\zeta}}
\newcommand{\uiota}{\underline{\iota}}
\newcommand{\ugamma}{\underline{\gamma}}

\newcommand{\simeqd}{\mathrel{\rotatebox[origin=c]{-90}{$\xrightarrow{\sim}$}}}

\newcommand{\seq}{\boldsymbol{\mu}}

\newcommand{\frg}{\mathfrak{g}}

\newcommand{\op}{{\mathrm{op}}}

\newcommand{\Hf}{\frac{1}{2}}
\newcommand{\hf}{\frac{1}{2}}

\newcommand{\ufv}{{\operatorname{uf}}}
\newcommand{\fv}{{\mathrm f}}

\newcommand{\clAlg}{\mathcal{A}}
\newcommand{\bClAlg}{\overline{\clAlg}}
\newcommand{\upClAlg}{\mathcal{U}}
\newcommand{\bUpClAlg}{\overline{\upClAlg}}
\newcommand{\alg}{\mathsf{A}}

\newcommand{\ud}{\underline{d}}
\newcommand{\uk}{\underline{k}}

\renewcommand{\tw}{\operatorname{tw}}

\newcommand{\frz}{\mathfrak{f}}

\ifdefined\supp
\else
\newcommand{\supp}{\operatorname{supp}}
\fi

\ifdefined\tw
\else
	\newcommand{\tw}{{\opname tw}}
\fi

\ifdefined\G
    \renewcommand{\G}{{\mathbb G}}
\else
   \newcommand{\G}{{\mathbb G}}
\fi

\ifdefined\C
    \renewcommand{\C}{{\mathbb C}}
\else
   \newcommand{\C}{{\mathbb C}}
\fi

\usepackage{comment}

\author{Fan Qin}
\address{School of Mathematical Sciences \\ Beijing Normal University \\ China}
\email{qin.fan.math@gmail.com}

\author{Milen Yakimov}
\address{Department of Mathematics \\ Northeastern University \\ 360 Huntington Ave, Boston, MA 02115 \\ USA 
and International Center for Mathematical Sciences, Institute of Mathematics and Informatics \\
Bulgarian Academy of Sciences \\ 
Acad. G. Bonchev Str., Bl. 8 \\
Sofia 1113, Bulgaria
}
\email{m.yakimov@northeastern.edu}


\title[Cluster structures for simple algebraic groups]{
Partially compactified 
quantum cluster structures on simple algebraic groups and the full Berenstein--Zelevinsky conjecture
}

\begin{document}
\begin{abstract}
The construction of partially compactified cluster algebras on coordinate rings is handled by using codimension 2 arguments on cluster covers. An analog of this in the quantum situation is highly desirable but has not been found yet. In this paper, we present a general method for the construction of partially compactified quantum cluster algebra structures on quantized coordinate rings from that of quantum cluster algebra structures on localizations.  As an application, we construct a partially compactified quantum cluster algebra structure on the quantized coordinate ring of every connected, simply connected complex simple algebraic group. Along the way, we settle in full the Berenstein--Zelevinsky conjecture that all quantum double Bruhat cells have quantum cluster algebra structures associated to seeds indexed by arbitrary signed words, and prove that all such seeds are linked to each other by mutations. 
\end{abstract}

\maketitle
\setcounter{tocdepth}{2} 
	
\section{Introduction}	\label{sec:intro}
\subsection{Partially compactified quantum cluster structures on quantized coordinate rings}
\label{sec:int1}
Cluster algebras were defined by Fomin and Zelevinsky \cite{fomin2002cluster} in 2001 and quickly started playing a fundamental role in many areas of mathematics and mathematical physics. Their quantum counterparts were defined by Berenstein and Zelevinsky in \cite{BerensteinZelevinsky05}. One of the major motivations that was put forward by Fomin and Zelevinsky (already in the introduction of \cite{fomin2002cluster}) was a program for studying canonical bases: (Step 1) prove that the coordinate rings of important varieties and their quantizations possess cluster structures and (Step 2) study the canonical bases of these (quantized) coordinate rings via cluster algebra techniques. This program has played a major role in the subject and canonical bases of (quantum) cluster algebras have been developed in broad generality \cite{fock2006moduli,GeissLeclercSchroeer10b,musiker2013bases,gross2018canonical,davison2019strong,qin2017triangular,qin2019bases}. 

For step 1, the construction of (quantum) cluster algebra structures on (quantized) coordinate rings, one has to use partially compactified (quantum) upper cluster algebras, meaning cluster algebras in which the frozen variables are not inverted. (The term was coined in \cite{gross2018canonical}.) This is so because (quantum) upper cluster algebras, whose frozen variables are inverted are often associated to the (quantized) coordinate rings of Zariski open subsets of the varieties in question. 

On the classical level the construction of an upper cluster algebra structure or a partially compactified upper cluster algebra structure on the coordinate ring of a variety $X$ is often carried out by means of coverings of $X$ by tori and products of tori and affine spaces related by mutations, and then by invoking codimension 2 arguments for the complements of their unions, see \cite{BerensteinFominZelevinsky05,GSV-plethora} for the two situations, respectively.

It is desirable to have quantum analogs of this and other geometric methods for constructing cluster algebra structures. Although there is a widely used notion of Gelfand--Kirillov dimension in noncommutative algebra, the search for quantum analogs of these methods has been unsuccessful. However, there are many constructions of partially compactified quantum cluster algebra structures on various iterated skew polynomial extensions, e.g. \cite{GeissLeclercSchroeer11,GY13,Kang2018,goodearl2016berenstein,bittmann,GY21}. We can localize at the frozen variables to cook up quantum cluster algebra structures on the quantized coordinate rings of Zariski open subsets of important varieties. So now the problem becomes:
\medskip

\noindent
{\bf{General Problem.}} {\em{How do we construct a partially compactified upper quantum cluster algebra structure on a noncommutative algebra $A$ from the knowledge of a quantum cluster algebra structure on a localization of $A$?}}
\medskip

In particular, we have the following key subproblem which has been much discussed in the literature but was still open (away from type $A$):
\medskip

\noindent
{\bf{Specific Problem.}} {\em{How do we construct a partially compactified upper quantum cluster algebra structure on the quantized coordinate ring $R_q[G]$ of a complex simple algebraic group $G$ from the known quantum cluster algebra structure on the quantized coordinate ring $R_q[G^{w_0,w_0}]$
of the open double Bruhat cell of $G$?
}}

In the special case of the quantized coordinate ring of a Grassmannian the General Problem was resolved by Grabowski and Launois in \cite{GrabLaun} by invoking arguments with $\Z$-gradings of cluster algebras.

In this paper we present a general method for solving the General Problem and a resolution of the Specific Problem.

\subsection{Statements of main results} Along the way of resolving the Specific Problem we prove the full Berenstein--Zelevinsky conjecture \cite{BerensteinZelevinsky05} about quantum 
cluster algebra structures on the quantized coordinate rings of double Bruhat cells. The Berenstein--Zelevinsky seeds in question will also be used for the construction of a partially compactified upper quantum cluster algebra structure on the quantized rings $R_q[G]$. We briefly recall those seeds next. 
For two integers $m \leq n$, denote
\[
[m,n]= \{m, m+1, \ldots, n\}. 
\]

Fix a (connected, simply connected) complex simple algebraic group $G$ and denote by $R_q[G]$ the quantized coordinate ring of $G$ over $\C(q^\hf)$. 
Denote by $r$ the rank of $G$ and by $W$ its Weyl group. Let $l : W \to \N$ be the length function on $W$. 

Given a pair of Weyl group elements $(u,w)$, consider the {\em{double Bruhat cell}} 
\[
G^{u,w} := B^+ u B^+ \cap B^- w B^- \subset G,
\]
where $B^\pm$ are a pair of opposite Borel subgroups of $G$. Denote by $R_q[G^{u,w}]$ the quantized coordinate ring of $G^{u,w}$ over $\C(q^{\Hf})$.
A {\em{signed word}} $\ubi$ for $(u,w)$ is a word 
\begin{equation}
\label{eq:signed-word}
\bi_{1},\bi_{2}, \ldots, \bi_l
\end{equation}
in the alphabet $[-r,-1]\sqcup [1,r]$ such that $\ubi$ is a reduced word for the element $(u, w) \in W \times W$, where the simple reflections of the first copy of $W$ are enumerated by $[-r,-1]$ and those for the second copy by $[1,r]$. For $k \in [1,l]$, let 
$(u_{\leq k}, w_{\leq k}) \in W \times W$ be the element that corresponds to the product of the simple reflections for the subword $\bi_{1}, \ldots, \bi_k$, see \cite[Section 10.1]{BerensteinFominZelevinsky05}. 
Denote $\ddI:= [-r, -1]\sqcup [1,l]$. Choose and fix any permutation $(\bi_{-r},\ldots,\bi_{-2},\bi_{-1})$ of $[1,r]$. For $k \in I$, set 
\[
\gamma_k = 
\begin{cases}
\vpi_{\bi_{k}}, & \mbox{if} \, \, k \in [-r,-1]
\\
u_{\leq k} \vpi_{|\bi_k|}, & \mbox{if} \, \, k \in [1,l]
\end{cases}
\quad \quad \mbox{and} \quad \quad 
\delta_k = 
\begin{cases}
w^{-1} \vpi_{\bi_{k}}, & \mbox{if} \, \, k \in [-r,-1]
\\
w^{-1} w_{\leq k} \vpi_{|\bi_k|}, & \mbox{if} \, \, k \in [1,l],
\end{cases}
\]
where $\vpi_s$ denote the associated fundamental weights for $s\in[1,r]$, see Section \ref{subsec:Lie-alg}.

Berenstein and Zelevinsky \cite{BerensteinZelevinsky05} associated a (candidate) quantum seed $\sd^{\BZ}=\sd^{\BZ}(\ubi)$ of $R_q[G^{u,w}]$ to every signed word $\ubi$ for $(u,w)$. (For technical reasons, it will be denoted by the opposite seed $\dsd(\ubi)^\op$ in our conventions in Section \ref{sec:ca-sw}). The cluster variables in it are the generalized quantum minors (see \eqref{eq:generalized-minor})
\[
x_k := \Delta_{\gamma_k, \delta_k}, \quad k \in \ddI.
\]
The $\Lambda$-matrix of $\sd^{\BZ}(\ubi)$ is the skew-symmetric integer 
$\ddI \times \ddI$-matrix with entries
\[
\lambda_{kj} = \langle \gamma_k, \gamma_j \rangle - 
\langle \delta_k, \delta_j \rangle \quad 
\mbox{for} \quad k > j. 
\]
The exchange matrix of $\sd^{\BZ}(\ubi)$ equals minus the matrix in  \eqref{eq:signed-word-B-matrix}. The frozen variables of $\sd^{\BZ}(\ubi)$ are the ones indexed by $[-r,-1]$ and the maximal elements of all subsets $\{ j \in [1,l] \mid |\bi_j| = k \} \subset [1,l]$ for a fixed value of $k \in [1,r]$. Berenstein and Zelevinsky \cite{BerensteinZelevinsky05} proved that $\sd^{\BZ}(\ubi)$ is indeed a quantum seed of the skew field of fractions of the associated quantum torus embedded in the skew field of fractions of $R_q[G^{u,w}]$ but it remained to be proved that the two division algebras coincide. 

We will say that the signed word \eqref{eq:signed-word} is {\em{unshuffled}} if 
$\bi_1, \ldots, \bi_{l(u)}$ and $-\bi_{l(u)+1}, \ldots, -\bi_l$ are reduced words for $u$ and $w$, respectively, in the alphabet $[1,r]$. 

For a quantum seed $\sd$, $\clAlg(\sd)$ and $\upClAlg(\sd)$ will denote the quantum cluster and upper cluster algebras associated to $\sd$ (all frozen variables inverted) and 
$\bClAlg(\sd)$ and $\bUpClAlg(\sd)$ will denote 
the {\em{partially compactified}} quantum cluster and upper cluster algebras associated to $\sd$ (all frozen variables not inverted).
We refer the reader to Section \ref{sec:CA} for background on cluster algebras.

Berenstein and Zelevinsky conjectured that 
\[
R_q[G^{u,w}] = \upClAlg(\sd^{\BZ}(\ubi)) 
\]
for all signed words $\ubi$ for $(u,w)$, see \cite[Conjecture 10.10]{BerensteinZelevinsky05}. For unshuffled signed words, this was proved in \cite[Main Theorem]{goodearl2016berenstein} in a stronger form than originally conjectured that involved an equality with the associated quantum cluster algebras:
\[
R_q[G^{u,w}] = \upClAlg(\sd^{\BZ}(\ubi)) = \clAlg(\sd^{\BZ}(\ubi)).
\]
However, the full Berenstein--Zelevinsky conjecture for arbitrary signed words remained open. Furthermore, the relationship between different Berenstein--Zelevinsky seeds remained unclear as well. Our first main result settles both problems:
\medskip

\noindent
{\bf{Theorem A.}} {\em{Let $G$ be a (connected, simply connected) complex simple algebraic group and $(u,w)$ a pair of Weyl group elements. Up to permutations of vertices, any two Berenstein--Zelevinsky quantum seeds for $R_q[G^{u,w}]$ associated to signed words for $(u,w)$ can be obtained from each other by successive mutations. All such seeds $\sd^{\BZ}(\ubi)$ are quantum seeds of $R_q[G^{u,w}]$ and 
\[
R_q[G^{u,w}] = \upClAlg(\sd^{\BZ}(\ubi)) = \clAlg(\sd^{\BZ}(\ubi)).
\]
}}
\medskip

The theorem is proved in Section \ref{sec:dBruhat}. 
Our proof is based on combining 
\begin{enumerate}
\item the special case of the conjecture proved in \cite{goodearl2016berenstein} and the cluster structure of CGL extensions \cite{GY13},
\item the structure results for quantum double Bruhat cells developed in \cite{joseph-paper,HLT,yakimov2014spectra},
\item the theory of common triangular bases developed in \cite{qin2017triangular,qin2020analog, qin2020dual,qin2023analogs,qin2024infinite},
\item the mutation results of Shen--Weng \cite{shen2021cluster} for classical seeds on double Bott--Samelson cells, and 
\item linking the Berenstein--Zelevinsky seeds associated to shuffled signed words to those for unshuffled via mutations and permutations using a result of Gei\ss, Leclerc and Schr\"oer \cite{GeissLeclercSchroeer11}.
\end{enumerate} 
\medskip

For the longest element $w_0$ of $W$, the quantized coordinate ring $R_q[G^{w_0,w_0}]$ is a localization of $R_q[G]$. The Berenstein--Zelevinsky seeds for $R_q[G^{w_0,w_0}]$ are seeds of $R_q[G]$ as well. Using Theorem A we prove the second main result of the paper
which settles the Specific Problem in Section \ref{sec:int1}. 
\medskip

\noindent 
{\bf{Theorem B.}} {\em{Let $G$ be a (connected, simply connected) complex simple algebraic group.  

The quantized coordinate ring of $G$ coincides with the partially compactified upper quantum cluster algebra associated to every 
Berenstein--Zelevinsky seed for any signed word $\ubi$ of $(w_0,w_0)$:
\[
R_q[G] = \bUpClAlg(\sd^{\BZ}(\ubi)).
\]
}}

We note that by Theorem A, the seeds $\sd^{\BZ}(\ubi)$ in Theorem B can be reached from each other by successive mutations. This result is proved in Theorem \ref{thm:quantized_coordinate_ring} (plus the mutation result in Theorem A). 

In Theorem \ref{thm:classical-cluster-G} in Appendix B, we prove a classical version of this result based on codimension $2$ arguments, namely that the coordinate ring $\C[G]$ of $G$ is isomorphic to the 
classical partially compactified upper cluster algebra associated to 
$\bUpClAlg(\sd^{\BZ})$. While completing this paper, the authors learned that Oya \cite{Oya} independently obtained a classical version of this result when $G$ is not of type $F_4$ using a different method and he also proved that $\bClAlg= \bUpClAlg$ in this case. The $A$ case was obtained earlier by Fomin--Williams--Zelevinsky
\cite[Section 6.6]{fomin2020introduction}. 
\medskip

The proof of Theorem B uses the following general method for constructing partially compactified upper quantum cluster algebra structures on quantized coordinate rings. This is the technique that we propose to address the General Problem in Section \ref{sec:int1}. 
\medskip

\noindent
{\bf{Theorem C.}} {\em{Assume that $A$ is a (noncommutative) algebra without zero divisors and $x_j \in A$ (for $j \in I_\fv$) are finitely many prime elements of $A$ (cf. Definition \ref{def:prime}) such that 
\begin{equation}
\label{eq:loc}
A[x_j^{-1}, j \in I_\fv] \simeq \upClAlg
\end{equation}
for some upper quantum cluster algebra $\upClAlg$ such that}}
\begin{enumerate} 
\item 
{\em{$A \subset \bUpClAlg$ for the partially compactified upper quantum cluster algebra $\bUpClAlg$ associated to $\upClAlg$, and}}
\item
{\em{the initial cluster variables of $\upClAlg$ belong to $A$ 
and the frozen ones are precisely $\{x_j \mid j \in I_\fv\}$.}}
\end{enumerate}
{\em{Then $A = \bUpClAlg$.
}}
\medskip

\noindent
{\bf{Remark.}} 
\begin{enumerate}
\item {\em{In concrete applications of Theorem C, we always have an inclusion of the form $A \subset \upClAlg$ that we start with. 
In particular, we have for free that $A$ has no zero divisors (i.e., it is a domain) because this is true for every upper quantum cluster algebra.}}
\item {\em{The inclusion $A \subset \bUpClAlg$ can then be verified effectively by using the
following characterization of $\bUpClAlg$. Denote by 
$\{x_j \mid j \in I \}$ the initial cluster of $\upClAlg$. 
Let $\bLP$ be the subalgebra of the quantum torus in the generators $x_j$ that is spanned by monomials for which the powers of the frozen variables $\{x_j \mid j \in I_\fv\}$ are nonnegative (no condition on the powers of the unfrozen variables). Then, 
\[
\bUpClAlg = \upClAlg \cap \bLP,
\]
cf. Corollary \ref{cor:Uint}. 
In other words, checking if an element of $\upClAlg$ belongs to 
$\bUpClAlg$ only amounts to checking a condition for one seed of 
the quantum cluster algebra.}} 
\end{enumerate}

\subsection*{Conventions} If $I=I_1\sqcup I_2$, we will view $\Z^{I_1}$, $\Z^{I_2}$ as subsets of $\Z^I$ via the natural decomposition $\Z^I:=\Z^{I_1}\oplus \Z^{I_2}$. All vectors are assumed to be column vectors. If $B=(b_{ij})_{i,j\in I}$ is an $I\times I$ skew-symmetric integer matrix, its adjacency quiver $Q$ is defined such that the set of vertices is $I$ and the number of arrows from $i$ to $j$ is $\max(0,b_{ij})$. We allow $A=B$ when we denote $A\subset B$.
\subsection*{Acknowledgements} We are grateful to Geoffrey Janssens and Malfeit Wannes for very helpful comments on the impact of the choice of frozen variables on the equivalence $\clAlg=\upClAlg$ 
$\Leftrightarrow$ $\bClAlg=\bUpClAlg$ in Proposition \ref{prop:AeqU} and an example in \cite{BMS}. We are indebted to the anonymous referees for their careful and insightful comments, which greatly helped us to improve the quality of this manuscript.

The research of Fan Qin has been supported by the National Natural Science Foundation of China (Grant No. 12422102) and that of Milen Yakimov by the Bulgarian Science Fund grant KP-06-N92/5 and the Ministry of Education and Science grant DO1-239/10.12.2024, the Simons Foundation grant SFI-MPS-T-Institutes-00007697 and the USA National Science Foundation grant DMS–2200762.
\section{Classical and quantum cluster algebras}
\label{sec:cl}
In this section we gather background material on cluster algebras.

\subsection{Basics of cluster algebras}
\label{sec:CA}

\subsubsection*{Seeds}

We choose and fix a finite set $I$ endowed with a partition $I=I_\ufv\sqcup I_\fv$. The elements of $I_\ufv$ and $I_\fv$ will be called unfrozen vertices and frozen vertices, respectively. Choose strictly positive relatively prime integers $\symm_i$, $i\in I$. 

Let $(b_{ij})_{i,j\in I}$ denote any matrix over $\frac{\Z}{2}$ such that $\symm_i b_{ij} =-\symm_j b_{ji}$, $\forall i, j \in I$. Assume that the submatrix $\tB:=(b_{ik})_{i\in I,k\in I_\ufv}$ has integer entries. Let $(x_i)_{i\in I}$ be indeterminates. 
\begin{definition}[\cite{fomin2002cluster}]
    The collection $\sd:=(\tB, (x_i)_{i\in I})$ is called a classical seed or simply a seed, $x_i$ are called cluster variables, and $\tB$ the exchange matrix of $\sd$. 
\end{definition}

Assume that there exists a skew-symmetric $\Z$-matrix $\Lambda=(\Lambda_{ij})_{i,j\in I}$ and strictly positive integers $\ddiag_k$ such that $\sum_{j\in I} \Lambda_{ij}b_{jk}=-\delta_{i,k}\ddiag_k$ for $k\in I_\ufv$. We say that the matrices $\Lambda$ and $\tB$ are compatible. It follows that in this case $\tB$ is of full rank.
\begin{definition}[\cite{BerensteinZelevinsky05}]
    The pair $(\tB,\Lambda)$ is called a compatible pair and $\sd:=(\tB,\Lambda,(x_i)_{i\in I})$ is called a quantum seed or a quantization of the classical seed $(\tB, (x_i)_{i\in I})$. $\Lambda$ is called the $\Lambda$-matrix or the quantization matrix of $\sd$.
\end{definition}
Denote $M:=\Z^I$ and its standard basis vectors by $f_i$, $i\in I$. The lattice associated to the seed $\sd$ will be denoted by $M(\sd)$. 
The matrix $\Lambda$ determines a skew-symmetric bilinear form $\lambda$ on $M(\sd)$ given by 
\begin{align*}
    \lambda(g,h):=g^T \Lambda h,
\end{align*}
where $(\ )^T$ denotes the transpose.

Unless otherwise specified, we choose the base ring $\kk$ to be $\Z$ for classical seeds and $\Z[q^{\pm\Hf}]$ for quantum seeds, where $q^{\Hf}$ is a formal quantum parameter. The torus algebra associated with a seed $\sd$ is the Laurent polynomial ring 
\begin{align*}
    \LP(\sd):=\kk[x_i^{\pm 1}]_{i\in I}.
\end{align*}
Its commutative product will be denoted by $\cdot$. 
The algebra $\LP(\sd)$ contains the mixed torus-affine algebra
\begin{align*}
    \bLP(\sd):=\kk[x_j]_{j\in I_\fv}[x_k^{\pm 1}]_{k\in I_\ufv}.
\end{align*}

Define $N_{\ufv}:=\Z^{I_\ufv}$. For $m=(m_i)_{i\in I}\in M(\sd)$ and $n=(n_k)_{k\in I_\ufv}\in N_{\ufv}$, denote the elements
\[
x^m:=\prod_{i\in I} x_i^{m_i}, 
\quad 
y_k := x^{\sum_i b_{ik} f_i},
\quad \mbox{and} \quad  y^n:=\prod_{k\in I_\ufv} y_k^{n_k} \in \LP(\sd)
\]
with the commutative product on Laurent monomials. 

\begin{definition}
$x^m$ is called a cluster monomial of $\sd$ if $m\in \N^I$ and a localized cluster monomial if $m\in \N^{I_\ufv}\oplus \Z^{I_\fv}$. 

$x_j$, $j\in I_\fv$, are called frozen variables. A frozen factor is a Laurent monomial $x^m$ such that $m\in \Z^{I_\fv}$. We denote $\frGp:=\{x^m|m\in \Z^{I_\fv}\}$.
\end{definition}

Next, assume that $\sd$ is a quantum seed. We introduce the $\Z$-linear automorphism $\overline{(\ )}$ on $\LP(\sd)$, called the bar involution, given by  
\begin{align*}
    \overline{q^\alpha x^m}=q^{-\alpha}x^m.
\end{align*}
We also introduce the following $q$-twisted product on $\LP(\sd)$ and $\bLP(\sd)$:
\begin{align*}
    x^g * x^h:=q^{\Hf\lambda(g,h)}x^{g+h}.
\end{align*}
The product $*$ will be our default multiplication on $\LP(\sd)$ and $\bLP(\sd)$, whose symbol will often be omitted for simplicity. Note that the classical product on $\bLP(\sd)$ is given by the same formula with $q^{\Hf}$ replaced by $1$. We use $\cF$ to denote the skew-field of fractions of $\LP(\sd)$. 

Let $\sigma$ denote any permutation of $I$ such that $\sigma I_\ufv=I_\ufv$, and thus $\sigma I_\fv = I_\fv$. For any matrix $Z=(z_{ij})$ whose rows and columns are indexed by $I$, $I_\ufv$, or $I_\fv$, we introduce the matrix $\sigma Z$ such that $(\sigma Z)_{\sigma i,\sigma j}=z_{ij}$. We then define the permuted seed $\sigma \sd$ by setting $\tB(\sigma \sd):=\sigma \tB$, $x_{\sigma i}(\sigma \sd):=x_i$, and, in the quantum case, $\Lambda(\sigma \sd):=\sigma \Lambda$. If $\sigma$ is a permutation on $I_\ufv$, we extend its action to $I$ by acting on $I_\fv$ trivially. If $\sd'=\sigma\sd$, we will often identify $\sd'$ with $\sd$ via $\sigma$ and simply denote $\sd'=\sd$.

\subsubsection*{Mutations}
Consider a (classical or quantum) seed $\sd$. Following \cite{fomin2002cluster}\cite{BerensteinFominZelevinsky05}, for any unfrozen vertex $k$, a mutation operation $\mu_k$ is defined, which produces a new seed $\sd':=\mu_k \sd:=(\tB',(x_i')_{i\in I})$. If $\sd$ is a quantum seed, we obtain a quantum seed $\sd'$ endowed with a $\Lambda$-matrix $\Lambda'$.

More precisely, denote $[ b]_+ := \max(0, b)$ and set $\tB':=(b'_{ij})_{i\in I,j\in I_\ufv}$, where the entries $b'_{ij}$ are given by
\begin{align*}
    b'_{ij} & :=\begin{cases}
    b_{ij}+[b_{ik}]_+[b_{kj}]_+ - [-b_{ik}]_+[-b_{kj}]_+, & i,j\neq k\\
    -b_{ij}, & i=k \text{ or }j=k.
    \end{cases}
\end{align*}
Let us define the piecewise linear map $\phi_{\sd',\sd}:M(\sd)\rightarrow M(\sd')$ sending $m=(m_i)_{i\in I}\in \Z^I$ to $m'=(m'_i)_{i\in I}\in \Z^I$ such that
\begin{align*}
    m'_k & = -m_k,\\
    m'_i & = m_i + [b_{ik}]_+[m_k]_+ - [-b_{ik}]_+[-m_k]_+, \forall i\neq k.
\end{align*}
Note that $\phi_{\sd,\sd'}$ and $\phi_{\sd',\sd}$ are inverse to each other.
Denote the standard basis of $M(\sd')$ by $\{f'_i\mid i\in I\}$. The $\Lambda$-matrix of the seed $\sd'$ is the matrix $\Lambda':=(\Lambda'_{ij})_{i,j\in I}$ with entries
\begin{align*}
    \Lambda'_{ij} & =\lambda(\phi_{\sd,\sd'}f'_i,\phi_{\sd,\sd'}f'_j).
\end{align*}
The skew-field of fractions of $\LP(\sd)$ will be denoted by $\cF'$. We have the $\kk$-algebra isomorphism $\mu^*_{\sd',\sd}:=\mu_k^*:=\cF'\rightarrow \cF$, such that 
\begin{align*}
    \mu_k^*(x_i')&=x_i,\quad \forall i\neq k,\\
    \mu_k^*(x'_k)&=x^{-f_k + \sum_j [-b_{jk}]_+f_j}+x^{-f_k + \sum_i [b_{ik}]_+f_i}.
\end{align*}
Note that we have $\mu_k^*(x'_k)=x^{\phi_{\sd,\sd'}f_k}\cdot (1+y_k)$. We will identify $\cF'$ and $\cF$ via the isomorphism $\mu_k^*$ and thus view $x_i'$ as elements of $\cF$. We will often omit the symbol $\mu_k^*$ for simplicity.

It is well-known that $\mu_k$ is an involution, i.e., $\mu_k( \mu_k \sd)=\sd$. If $\uk:=(k_1,\ldots,k_r)$ is a sequence of unfrozen vertices, denote the mutation sequence $\seq_{\uk}:=\mu_{k_r}\cdots \mu_{k_1}$ and the seeds $\sd_t:=\mu_{k_t}\cdots \mu_{k_2}\mu_{k_1}\sd$ for $0\leq t\leq r$, where $\sd_{0}:=\sd$. Define the mutation isomorphism 
\[
\mu_{\sd_r,\sd}^*:=\seq_{\uk}^*=\mu_{k_1}^*\cdots \mu_{k_r}^*:\cF(\sd_r)\rightarrow \cF(\sd)
\]
of skew-fields. By \cite{gross2013birational}, computing the tropicalization of $\seq_{\uk}^*$, we obtain the composition 
\[
\phi_{\sd_r,\sd}=\phi_{\sd_r,\sd_{r-1}}\cdots \phi_{\sd_2,\sd_1}\phi_{\sd_1,\sd_0}:M(\sd)\rightarrow M(\sd_r),
\]
which depends only on the seeds $\sd$ and $\sd_r$, 
but not on the sequence $k_1, \ldots, k_r$.

For any $z\in \LP(\sd)$ and $j\in I$, let $\nu_j^\sd(z)$ denote its order of vanishing at $x_j=0$.
\begin{lem}[{\cite[Lemma 2.12]{qin2023analogs}}]\label{lem:vanishing-order-inv}
    For any $k\in I_\ufv$ and $z\in \LP(\sd)\cap \LP(\mu_k \sd)$, we have $\nu_j^\sd(z)=\nu_j^{\mu_k \sd}(z)$ whenever $j\neq k$.
\end{lem}

\subsubsection*{Cluster algebras}
In what follows we work simultaneously with classical and quantum cluster algebras, where the first setting uses the commutative product of $\LP(\sd)$ and the second the noncommutative one. An extra $q$ will not be placed in the notation for quantum algebras, instead the classical and quantum picture will be distinguished by the initial classical or quantum seed.  

Let $\Delta^+:=\Delta^+_\sd$ denote the set of seeds obtained from an initial (classical or quantum seed) 
$\sd$ by finitely many mutations. Following \cite{fomin2002cluster}, \cite{BerensteinFominZelevinsky05}, and \cite{gekhtman2012cluster}, we define cluster algebras as the following.
\begin{definition}
    The (partially compactified) cluster algebra is defined to be 
    \[
    \bClAlg:=\bClAlg(\sd):=\kk[x_i(\sd')]_{\sd'\in \Delta^+_\sd}.
    \]
    The (partially compactified) upper cluster algebra is defined to be 
    \[
    \bUpClAlg:=\bUpClAlg(\sd):=\bigcap_{\sd'\in \Delta^+_\sd}\bLP(\sd').
    \]
    Their localization $\clAlg:=\clAlg(\sd):=\bClAlg(\sd)[x_j^{-1}]_{j\in I_\fv}$ and $\upClAlg:=\upClAlg(\sd):=\bUpClAlg(\sd)[x_j^{-1}]_{j\in I_\fv}$ are called the (localized) cluster algebra and the (localized) upper cluster algebra, respectively.
\end{definition}
By a cluster algebra $\alg$, we might mean $\bClAlg$, $\bUpClAlg$, $\clAlg$, or $\upClAlg$.

By Lemma \ref{lem:vanishing-order-inv}, for $j\in I_\fv$, the order of vanishing $\nu_j^{\sd'}$ is independent of the choice of seed $\sd'\in \Delta^+$. So we can simply denote it by $\nu_j$. 
Since $\bLP(\sd)=\{z\in \LP(\sd)|\nu_j(z)\geq 0,\ \forall j\in I_\fv \}$, we have that
\begin{align}\label{eq:order-bUpClAlg}
\bUpClAlg=\{z\in \upClAlg|\nu_j(z)\geq 0,\ \forall j\in I_\fv \}=\upClAlg\cap \bLP(\sd).
\end{align}

The seed opposite to $\sd$ is defined as $\sd^\op:=(-\tB,(x_i)_{i\in I})$, and $\Lambda(\sd^\op):=-\Lambda(\sd)$ when $\sd$ is a quantum seed. We have a $\kk$-algebra anti-isomorphism $\iota^\op:\clAlg(\sd)\rightarrow\clAlg(\sd^\op)$, sending $x_i(\sd)$ to $x_i(\sd^\op)$. We can verify $\iota^\op (x_i(\seq \sd))=x_i(\seq \sd^\op)$ for any mutation sequence $\seq$. Moreover, $\iota^\op$ identifies the bar-involutions, i.e., $\overline{\iota^\op z}=\iota^\op \overline{z}$, $\forall z\in \LP(\sd)$. See \cite{qin2020analog}.
\subsubsection*{Degrees and dominance orders}
Let us work with any (classical or quantum) seed $\sd$. Recall that $M(\sd)=\Z^I$ and $N_\ufv=\Z^I_{\ufv}$. Introduce $N^\oplus=\N^{I_\ufv}$ and $M^\oplus:=\tB N^\oplus$. 
When the seed $\sd$ needs to be emphasized we will use the notation $N^\oplus(\sd)$ and $M^\oplus(\sd)$.
\begin{definition}[{Dominance order \cite{qin2017triangular}\cite[Proof of Proposition 4.3]{cerulli2015caldero}}]
    For any $m,m'\in M(\sd)$, we say $m$ dominates $m'$, denoted $m'\preceq_\sd m$, if $m'\in m+M^\oplus$.
\end{definition}

\begin{definition}[\cite{qin2017triangular}]
An element $z\in\LP(\sd)$ is said to have degree $m$, for some $m\in M(\sd)$, if it takes the form $z=\sum_{m'\preceq_\sd m} b_{m'} x^{m'}$, $b_{m'}\in \kk$, such that $b_m\neq 0$. In this case, denote $\deg^{\sd} z:=m$. It is further said to be $m$-pointed if $b_m=1$. 

If $z$ has degree $m$ and $b_m\in q^{\Hf\Z}$, we define its normalization to be $[z]^\sd:=b_m^{-1}z$. 
\end{definition}
Note that any $m$-pointed element $z\in\LP(\sd)$ takes the form
    \begin{align*}
        z=x^m\cdot (1+\sum_{n\in N^\oplus(\sd)} c_n y^n),
    \end{align*} 
    where $c_n\in \kk$.

\begin{definition}
Let $\Theta$ denote a subset of $M(\sd)$. A subset $\cZ$ of $\LP(\sd)$ is said to be $\Theta$-pointed if it is of the form $\cZ=\{z_m|m\in \Theta\}$ with $z_m$ being $m$-pointed.
\end{definition}
Denote
\[
\mm:=q^{-\Hf}\Z[q^{-\Hf}].
\]
\begin{definition}
    Let $\cZ$ denote a $\Theta$-pointed set and $\sum b_m z_m$ a linear combination of its elements, where $b_m\in \kk$. If $\{m|b_m\neq 0 \}$ has a unique $\prec_\sd$-maximal degree, to be denoted $m^{(0)}$, and $b_{m^{(0)}}=1$, we say the element $\sum b_m z_m$ is $\prec_\sd$-unitriangular. 
    If, further, $b_{m}\in \mm$ whenever $m\neq m^{(0)}$, we say it is $(\prec_\sd,\mm)$-unitriangular.        
\end{definition}

We will often omit the symbol $\sd$ for simplicity. 

\begin{thm}[{\cite{FominZelevinsky07}\cite{Tran09}\cite{DerksenWeymanZelevinsky09}\cite{gross2018canonical}}]\label{thm:cluster-mono}
    All localized cluster monomials of $\clAlg(\sd)$ are pointed elements in $\LP(\sd')$ with distinct degrees or every seed $\sd'\in\Delta^+_{\sd}$.
    
    Moreover, assume that $\sd$ is a quantum seed. Let $z$ denote any quantum localized cluster monomial of $\clAlg(\sd)$. Then its semiclassical limit $z|_{q\mapsto 1}$ is the classical localized cluster monomial with the same degree.
\end{thm}

\subsubsection*{Common triangular basis }
Let us recall the common triangular bases in the sense of \cite{qin2017triangular}, which are analogs of the dual canonical bases for quantum cluster algebras, see \cite{qin2020dual}\cite{qin2023analogs}. 

\begin{definition}
    A seed $\sd$ is said to be injective-reachable if there exists a mutation sequence $\Sigma$ and a permutation $\sigma$ of $I$ such that $\deg x_{\sigma k}(\Sigma \sd)=-f_k$ for any $k\in I_\ufv$.   
\end{definition}
The mutation sequence $\Sigma$ is a green-to-red sequence in the sense of \cite{Keller11}.

Assume that $\sd$ is an injective-reachable seed, then so are all seeds in $\Delta^+_\sd$, see \cite{qin2017triangular}\cite{muller2015existence}. In this case, we denote $\sd[1]:=\Sigma \sd$. It is unique up to permutation, see \cite{FominZelevinsky07}\cite{gross2018canonical}.

Let $\alg$ denote $\clAlg$ or $\upClAlg$. 
\begin{definition}[\cite{qin2017triangular}]
    A triangular basis of a quantum cluster algebra $\alg$ with respect to the seed $\sd$ is a $\kk$-basis $\can$ such that:
    \begin{itemize}
        \item The elements of $\can$ are bar-invariant.
        \item $\can$ is $M(\sd)$-pointed: it takes the form $\can=\{\can_m|m\in M(\sd)\}$ with $\can_m$ being $m$-pointed.
        \item $\can$ contains the cluster monomials of $\sd$ and $\sd[1]$.
        \item For any $i\in I$ and $m\in M(\sd)$, we have the following $(\prec_\sd,\mm)$-unitriangular decomposition
        \begin{align}
            [x_i*\can_m]=\can_{m+f_i}+\sum_{m'\prec_\sd m+f_i} b_{m'} \can_{m'},\quad b_{m'}\in q^{-\Hf}\Z[q^{-\Hf}].
        \end{align}
    \end{itemize}
\end{definition}
By \cite{qin2017triangular}, a triangular basis with respect to $\sd$ is unique if it exists.

\begin{definition}[{\cite{qin2017triangular}\cite{qin2021cluster}}]
If a quantum cluster algebra $\alg$ has a $\kk$-basis $\can$ which is the triangular basis with respect to all seeds, $\can$ is called a common triangular basis. (It is unique if it exists.)
\end{definition}
By \cite{qin2020dual}, the triangular basis $\can$ with respect to $\sd$ is the common triangular basis if and only if it contains all cluster monomials of $\alg$. Moreover, if $\clAlg$ has a common triangular basis, it must be equal to the upper cluster algebra $\upClAlg$, see \cite{qin2019bases}.

\begin{definition}
    A frozen vertex $j\in I_\fv$ is said to be optimized in $\sd$ if $b_{jk}(\sd)\geq 0$ for all $k\in I_\ufv$. We say it can be optimized if it is optimized in some seed $\sd'=\seq \sd$ for some mutation sequence $\seq$.
\end{definition}

\begin{prop}[{\cite[Proposition 2.15]{qin2023analogs}}]\label{prop:optimized_bases}
Assume that the quantum upper cluster algebra $\upClAlg$ has a common triangular basis $\can$. If all frozen vertices can be optimized, then $\can\cap \bUpClAlg$ is a basis of $\bUpClAlg$, called the common triangular basis of $\bUpClAlg$.
\end{prop}

\begin{prop}[{\cite[Proposition 4.3.1]{qin2020analog}}]
If $\upClAlg(\sd)$ has a common triangular basis $\can$, then $\iota^\op \can$ is the common triangular basis for $\upClAlg(\sd^\op)$.
\end{prop}

\subsection{Quantum cluster algebras associated with signed words}
\label{sec:ca-sw}

\subsubsection*{Braids and signed words}
For any integers $j\leq k$, we use $[j,k]$ to denote the set of integers $s$ such that $j\leq s\leq k$. Let $r$ be a strictly positive integer and $C=(c_{ij})_{i,j\in [1,r]}$ denote a generalized Cartan matrix with symmetrizers $d_i$:
\begin{align*}
    c_{ii}&=2,\ \forall i\in [1,r]\\
    d_ic_{ij}&=d_j c_{ji},\ c_{ij}\in \Z_{\leq 0},\ \forall i\neq j.
\end{align*}
Let $\Br^+$ be the monoid of positive braids associated with $C$. It is generated by $\sigma_i$, $i\in [1,r]$, such that
\begin{align*}
    \sigma_i \sigma_j&=\sigma_j \sigma_i,\ \text{if }c_{ij}=0,\\
    \sigma_i \sigma_j \sigma_i&=\sigma_j \sigma_i \sigma_j,\ \text{if }c_{ij}c_{ji}=1,\\
    (\sigma_i \sigma_j)^2&=(\sigma_j \sigma_i)^2,\ \text{if }c_{ij}c_{ji}=2,\\
    (\sigma_i \sigma_j)^3&=(\sigma_j \sigma_i)^3,\ \text{if }c_{ij}c_{ji}=3.
\end{align*}
Let $e$ denote the identity element. The generalized braid group $\Br$ is generated by $\sigma_i^{\pm 1}$. The Weyl group $W$ is the quotient group of $\Br$ by the relations $\sigma_i^2=e$, $\forall i\in [1,r]$. Let $[\beta]$ denote the image of $\beta\in \Br$ in $W$. Denote $s_i:=[\sigma_i]$.

A word $\ueta=(\eta_1,\ldots,\eta_l)$ is a sequence formed from letters in $[1,r]$. It generates a positive braid $\beta_{\ueta}:=\sigma_{\eta_1}\cdots \sigma_{\eta_l}$ in $\Br^+$. Denote its length by $l(\beta_{\ueta}):=l(\ueta):=l$. A reduced expression of $w\in W$ is a factorization $w=s_{k_1}\ldots s_{k_r}$ of minimal length. Define the length of $w$ to be $l(w):=r$. We use $\ueta^\op$ to denote the opposite word associated with $\ueta$.

A signed word $\ubi$ is a sequence formed from letters in $[-r,-1]\sqcup [1,r]$. Let $\ueta$, $\uzeta$ denote any two words. Then any (riffle) shuffle of $\ueta$,$-\uzeta$ is a signed word. Conversely, any signed word is a shuffle of $\ueta$, $-\uzeta$ for some words $\ueta$, $\uzeta$.

Denote by $l(\ubi)$ the length of a signed word $\ubi$. For any $j\leq k\in [1,l(\ubi)]$, we let $\ubi_{[j,k]}$ denote the signed word $(\bi_{j},\bi_{j+1},\ldots,\bi_{k})$.

\subsubsection*{Seeds}
Following \cite{shen2021cluster}, we construct two seeds $\dsd=\dsd(\ubi)$ and $\rsd=\rsd(\ubi)$ for any given signed word $\ubi$. 

Set $l:=l(\ubi)$, $\dI:=\dI(\ubi):=[1,l]$, and $\ddI:=\ddI(\ubi):=[-r,-1]\sqcup \dI$. Denote $\bi_{k}:=|k|$ for $k\in [-r,-1]$. For any $k\in\ddI$, we introduce the following element in $\dI\sqcup\{+\infty\}$:
\begin{align*}
    k[1]:=\min(\{j \mid k<j\leq l,|\bi_k|=|\bi_j|\}\sqcup \{+\infty\}).
\end{align*}
We denote $k[0]=k$. We extend the definition of $k[d]$ recursively as follows: For any $d\in\N$, if $k[d]$ is already defined and is contained in $\ddI$, we define $k[d+1]:=(k[d])[1]$; similarly, if $k[-d]$ is already defined and contained in $\dI$, and there exists $k'\in\ddI$ such that $k[-d]=(k')[1]$, we define $k[-d-1]:=k'$.

When we have $j=k[d]\in \dI$ for $d\in \Z$ and $j[-1]\in[-r,-1]$, we denote $k^{\min}:=j$. Similarly, when we have $j=k[d]\in \dI$ for $d\in\Z$ and $j[1]=+\infty$, we denote $k^{\max}:=j$. We also define $|\bi_k|^{\max}:=j$ if $j=k^{\max}$.

Define the set of unfrozen vertices in $\dI(\ubi)$ and $\ddI(\ubi)$ to be $\dI_\ufv:=\dI_\ufv(\ubi):=\{k\in [1,l] \mid k[1]\leq l\}$. The sets of frozen vertices $\dI_\fv$ and $\ddI_\fv$ consist of the remaining vertices, respectively. For $k\in \ddI$, set $d_k:=d_{|\bi_k|}$, where $(d_i)_{i\in [1,r]}$ are the symmetrizers of the Cartan matrix introduced in Section \ref{subsec:Lie-alg}. 

Denote $\varepsilon_k=\sign(\bi_k)$. By extending the construction of exchange matrices of Berenstein--Fomin--Zelevinsky \cite{BerensteinFominZelevinsky05},
following \cite[Lemma 6.3]{qin2023analogs} and  \cite{shen2021cluster}, we introduce the $\ddI\times \dI_\ufv$-matrix $\ddB(\ubi)=(b_{jk})$ with entries
\begin{align}\label{eq:signed-word-B-matrix}
    b_{jk}=
    \begin{cases}
        \varepsilon_k & k=j[1] \\
        -\varepsilon_j & j=k[1]\\
        \varepsilon_k c_{|\bi_j|,|\bi_k|} & \varepsilon_{j[1]}=\varepsilon_k,\quad j<k<j[1]<k[1]\\
        \varepsilon_k c_{|\bi_j|,|\bi_k|} & \varepsilon_k=-\varepsilon_{k[1]},\quad j<k<k[1]<j[1]\\
        -\varepsilon_j c_{|\bi_j|,|\bi_k|} & \varepsilon_{k[1]}=\varepsilon_j,\quad k<j<k[1]<j[1]\\
        -\varepsilon_j c_{|\bi_j|,|\bi_k|} & \varepsilon_j=-\varepsilon_{j[1]},\quad k<j<j[1]<k[1]\\
        0 & \text{otherwise.}       
    \end{cases}.
\end{align}
Let $\dB$ denote its $\dI\times \dI_\ufv$-submatrix. These matrices have diagrammatic constructions, see \cite{shen2021cluster} or \cite[Section 6]{qin2023analogs}.

\begin{example}\label{eg:words-matrix}
Consider the case when $C$ is the Cartan matrix of type $A_2$ and 
$\ubi = (1,2,1,-1,-2,-1)$. In this case, we have $\dI=[1,6]$, $\ddI=[-2,-1]\sqcup[1,6]$, $(-2)[1]=2$, $(-1)[1]=1$, $1[1]=3$, $2[1]=5$, $3[1]=4$, $4[1]=6$, $5[1]=6[1]=+ \infty$, and consequently, $\ddI_\fv=\{-2,-1,5,6\}$, $\dI_\fv = \{5,6\}$. The matrix $\ddB(\ubi)$, whose rows are indexed, in order, by $(-2,-1,1,2,3,4,5,6)$
and whose columns are indexed, in order, by $(1,2,3,4)$ is given by
\[
\ddB(\ubi)=(b_{jk})=
\begin{pmatrix}
    -1&1&0&0\\
    1&0&0&0\\
    0&-1&1&0\\
    1&0&-1&1\\
    -1&1&0&-1\\
    0&-1&1&0\\
    0&1&0&-1\\
    0&0&0&1
\end{pmatrix}.
\]
We refer the reader to \cite[Example 2.9]{BerensteinFominZelevinsky05} or \cite[Example A.6]{qin2025introduction} for the quantum group nature of this setting.
\end{example}

\begin{definition}
    The seed $\dsd$ is the collection $(\ddB,(x_i)_{i\in \ddI})$. The seed $\rsd$ is the collection $(\dB,(x_i)_{i\in \dI})$. 
\end{definition}

For data associated with $\dsd(\ubi)$ or $\rsd(\ubi)$, we suppress the symbols $\dsd$ and $\rsd$ when no confusion can arise; when the choice of $\ubi$ is clear, we also suppress $\ubi$.

\begin{lem}[\cite{qin2023analogs}]\label{lem:full-rank}
The seeds $\dsd$ and $\rsd$ are of full rank; that is, their exchange matrices are of full rank.
\end{lem}
By Lemma \ref{lem:full-rank}, we can endow $\dsd$ and $\rsd$ with compatible $\Lambda$-matrices, see \cite{GekhtmanShapiroVainshtein03}\cite{GekhtmanShapiroVainshtein05}. 
Such compatible $\Lambda$-matrices are not unique; below we explain how to choose them in a consistent way.
By abuse of notation, the corresponding quantum seeds will be denoted by the same symbols; from the context it will be clear which one is being used.

\subsubsection*{Operations on signed words}

Denote $\rsd=\rsd(\ubi)$ and $\dsd=\dsd(\ubi)$. Following \cite[Section 2.3, Proposition 3.7]{shen2021cluster}, we introduce several operations on $\ubi$ to produce new signed words $\ubi'$ and new seeds $\rsd':=\rsd(\ubi')$, $\dsd'=\dsd(\ubi')$.

\begin{enumerate}
    \item (Left reflection) A left reflection on $\ubi$ produces a new signed word $(-\bi_{1},\bi_{2},\ldots,\bi_{l})$, denoted $\ubi'$. Note that we have $\rsd'=\rsd$.
    \item (Flips) Assume that $\varepsilon_k\neq \varepsilon_{k+1}$ for some $k\in [1,l]$. We can flip $\ubi$ to obtain a new signed word $\ubi'=(\ubi_{[1,k-1]},\bi_{k+1},\bi_{k},\ubi_{[k+2,l]})$. By \cite{shen2021cluster}, we are in one of the following two cases.
    \begin{itemize}
    \item Assume $|\bi_k| \neq |\bi_{k+1}|$, then $\dsd'=\dsd$ and $\rsd'=\rsd$.
    \item Assume $|\bi_k|=|\bi_{k+1}|$. Then $k\in I_\ufv(\ubi)$, and we have $\dsd'=\mu_k\dsd$, $\rsd'=\mu_k\rsd$, where no permutation of seeds are needed.
    \end{itemize}
    \item (Braid moves) Assume that $\ubi_{[j,k]}$ and $\ugamma$ are two words in $[1,r]$, which produce the same positive braid $\beta_{\ubi_{[j,k]}}=\beta_{\ugamma}$. Consider $\ubi'=(\ubi_{[1,j-1]},\ugamma,\ubi_{[k+1,l]})$. In this case, there exists a sequence of mutations $\seq_{\ubi',\ubi}$ involving the unfrozen variables in $U:=\{r\in [j,k]|r[1]\leq k\}$ and a permutation of $[j,k]$ satisfying $\sigma U=U$, such that $\dsd'=\sigma \seq_{\ubi',\ubi}\dsd$ and $\rsd'=\sigma \seq_{\ubi',\ubi}\rsd$. 
\end{enumerate}

As a special case, $\uiota=(\uzeta^\op,\ueta)$ can be obtained from $\ubi$ by a sequence of left reflections and flips, see \cite[Section 6]{qin2023analogs}. 

We will denote $\sigma \seq_{\ubi',\ubi}$ by $\seq^\sigma_{\ubi',\ubi}$, called a permutation sequence, where the superscript $\sigma$ should not be taken as a specified permutation -- it indicates that a permutation is involved.

Let $\ubi^{(j)}$, $j=1,2,3$, denote three signed words. Let $\sd^{(j)}$ denote $\rsd(\ubi^{(j)})$ or $\dsd(\ubi^{(j)})$. When $\sd^{(j)}=\rsd(\ubi^{(j)})$, we assume that $\ubi^{(j)}$ are connected by left reflections, flips, and braid moves; when $\sd^{(j)}=\dsd(\ubi^{(j)})$, we assume that $\ubi^{(j)}$ are connected by flips and braid moves. Let $\seq^\sigma_{\ubi^{(j')},\ubi^{(j)}}$ denote any chosen permutation-mutation sequences associated with these operations, such that $\sd^{(j')}=\seq^\sigma_{\ubi^{(j')},\ubi^{(j)}}\sd^{(j)}$. 

By \cite{shen2021cluster}, we have the following commutative diagram for classical seeds.
\begin{align}
    \begin{array}{ccc}
    \clAlg^{(3)} & \overset{(\seq_{\ubi^{(3)},\ubi^{(1)}}^{\sigma})^{*}}{\xrightarrow{\sim}} & \clAlg^{(1)}\\
    \simeqd(\seq_{\ubi^{(3)},\ubi^{(2)}}^{\sigma})^{*} &  & \parallel\\
    \clAlg^{(2)} & \overset{(\seq_{\ubi^{(2)},\ubi^{(1)}}^{\sigma})^{*}}{\xrightarrow{\sim}} & \clAlg^{(1)}
    \end{array}.\label{eq:connect-words}
    \end{align}

Next, endow $\sd^{(1)}$ with a compatible $\Lambda$-matrix. We then endow $\sd^{(j+1)}$ with a compatible $\Lambda$-matrix from that of $\sd^{(j)}$ via the permutation-mutation sequence $\seq^\sigma_{\ubi^{(j+1)},\ubi^{(j)}}$ for $j=1,2$. Then the quantum seed $\sd^{(3)}$ equals $\seq^\sigma_{\ubi^{(3)},\ubi^{(1)}}\sd^{(1)}$, and the above diagram is still commutative for quantum cluster algebras, see \cite[Lemma 3.2]{qin2024infinite}.

\subsubsection*{Properties}
We will consider quantum seeds $\dsd$ and $\rsd$ from now on.

\begin{lem}[\cite{qin2023analogs}]\label{lem:optimized}
The frozen vertices of $\dsd$ and $\rsd$  can be optimized.
\end{lem}
By Proposition \ref{prop:optimized_bases} and Lemma \ref{lem:optimized}, the common triangular basis of $\upClAlg(\dsd)$ restricts to that of $\bUpClAlg(\dsd)$, and the common triangular basis of $\upClAlg(\rsd)$ restricts to that of $\bUpClAlg(\rsd)$.

\begin{thm}[\cite{GY13,GY18,shen2021cluster,qin2023analogs}]\label{thm:A=U}
    We have $\clAlg(\rsd)=\upClAlg(\rsd)$ and $\clAlg(\dsd)=\upClAlg(\dsd)$.
    \end{thm}

\subsubsection*{Standard bases}

Let us take $\uiota=(\uzeta^\op,\ueta)$ as before. By \cite[Lemma 8.4]{qin2023analogs}, there is a green-to-red sequence $\Sigma$ from $\rsd:=\rsd(\uiota)$ to $\rsd[1]$, such that the cluster variables of the seeds appearing along $\Sigma$ take the following form:
\begin{itemize}
    \item The cluster variables appearing are denoted by $W_{[j,k]}$ for $1\leq j\leq k\leq l$.
    \item $\deg^{\rsd(\uiota)} W_{[j,k]}=f_k-f_{j[-1]}$, where we understand $f_{j[-1]}=0$ when $j=j^{\min}$.
\end{itemize} 
In fact a maximal green sequence with these properties exists by  
\cite[Theorem A]{Ymaxgreen}.
We will not give the original definition of $W_{[j,k]}$. It suffices to know that they are uniquely determined by the second property.

Note that the initial variables $x_k$ for $k\in [1,l]$ are denoted by $W_{[k^{\min},k]}$.
\begin{definition}
    The cluster variables $W_{[j,k]}$ are called the interval variables of $\upClAlg(\rsd(\uiota))$. $W_{k}:=W_{[k,k]}$ are called the fundamental variables. Denote $\gamma_{[j,k]}:=\deg W_{[j,k]}$ and $\gamma_k:=\deg W_k$. 
\end{definition}

\subsubsection*{Standard bases and Kazhdan--Lusztig bases}
Let $<_{\lex}$ denote the lexicographical order on $\Z^{[1,l]}$, where the coordinates are compared from left to right.
We also introduce the reverse lexicographical order such that, for $w,w'\in \N^{[1,l]}$, we have $w<_{\rev} w'$ if and only if $(w')^\op<_{\lex} w^\op$, where $(w_1,\ldots,w_l)^\op:=(w_l,\ldots,w_1)$. 

Note that $<_{\lex}$ and $<_{\rev}$ are bounded from below on $\N^{[1,l]}$. Recall that $\gamma_j=\deg W_j$.
\begin{prop}[{\cite[Lemma 8.8]{qin2023analogs}}]
    For any $w=(w_j),w'=(w'_j)\in \N^{[1,l]}$, if $\sum w'_j \gamma _j\prec_{\rsd} \sum w_j \gamma_j$, then $w'<_{\lex}w$ and $w'<_{\rev}w$.
\end{prop}

For any $\uc=(c_j)_{j\in [1,l]}\in \N^{[1,l]}$, we define the standard monomial $\std(\uc)$ as a normalized product of the fundamental variables in $\LP(\rsd)$:
\begin{align*}
    \std(\uc)= [W_1^{c_1}*W_2^{c_2} *\cdots* W_{l}^{c_{l}}]
\end{align*}
Note that $\std(0)=1$. Denote $\std:=\{ \std(\uc) \mid \uc\in \N^{[1,l]}\}$.
\begin{thm}[\cite{qin2023analogs}]\label{thm:bases_dBS}
(1) $\std$ is a $\kk$-basis of $\bUpClAlg(\rsd)$, called the standard basis.

(2) The standard basis satisfies the analog of the Levendorskii--Soibelman straightening law:
\begin{align*}
    W_k W_j-q^{\lambda(\deg W_k,\deg W_j)}W_j W_k\in \sum_{\uc\in \N^{[j+1,k-1]}}\kk \std(\uc).
\end{align*}

(3) We have $\bUpClAlg(\rsd)=\bClAlg(\rsd)$.

(4) Let $\can:=\{\can(\uc)|\uc\in \N^{[1,l]}\}$ denote the Kazhdan--Lusztig basis (KL basis for short) associated with $\std$ sorted by $<_{\lex}$, i.e., $\can(\uc)$ are bar-invariant and satisfy
\begin{align*}
    \can(\uc)=\std(\uc)+\sum_{\uc'<_{\lex}\uc}b_{\uc'}\std(\uc'),\quad b_{\uc'}\in \mm.
\end{align*}
Then $\can$ is the common triangular basis of $\bUpClAlg(\rsd)$. Moreover, the statement is still true if we replace $<_{\lex}$ by $<_{\rev}$. 
\end{thm}
Note that $\can(\uc)$ are $\sum c_j \gamma_j$-pointed.

\section{Quantum groups}
\label{sec:q-grps}
In this section we gather background material on quantum groups. 

\subsection{Semisimple Lie algebras}\label{subsec:Lie-alg}
Let $\frakg$ be a complex semisimple Lie algebra of rank $r$
with Weyl group $W$. Denote by $\Pi= \{ \alpha_1, \ldots, \alpha_r \}$,  $\{s_i\}$, 
$\{\alpha_i\spcheck 
\}$, $\{\varpi_i\}$ and $\{\varpi_i\spcheck\}$ the sets of its simple roots, simple reflections, simple coroots, fundamental weights, and fundamental coweights. Let $P$ and $Q$ be the weight and root lattices of $\frakg$, $P^+$ be the set of dominant integral weights of $\frakg$, and
$Q^+:= \sum_i \Z_{\geq 0} \alpha_i$. Let $\langle.,. \rangle$ be the invariant bilinear form on $\R \Pi$
normalized by $\langle \alpha_i, \alpha_i \rangle = 2$ for short roots $\alpha_i$. In this notation the Cartan matrix of $\frakg$ is  
\[
(c_{ij}) := ( \langle \alpha_i\spcheck, 
\alpha_j \rangle) \in M_r(\Z).
\]
Denote
$$
|| \gamma ||^2 
:= \langle \gamma, \gamma \rangle \quad 
\mbox{for} \quad \gamma \in \R \Pi.
$$
The symmetrizing integers of the Cartan matrix of $\frakg$
are $d_i := || \alpha_i ||^2/2$ . Let 
\[
\rho := \vpi_1 + \cdots + \vpi_r \quad \mbox{and}
\quad
\rho\spcheck := \vpi_1\spcheck + \ldots + \vpi_r\spcheck.
\]
The joint support $\cS(w_1, \ldots, w_m)$ of a set of Weyl group elements $w_1, \ldots, w_m \in W$ consists of those $i \in [1,r]$ such that $s_i$ appears in a reduced expression of some $w_j$.

\subsection{Quantum groups and quantum function algebras}
\label{sec:q-funct}
Throughout the rest of the section we work over $\K:= \C(q)$
or $\K:= \C(q^\Hf)$. Set $q_i=q^{d_i}$. The quantized universal enveloping algebra $U_q(\g)$ of $\g$ over $\KK$ has Chevalley generators $K_i^{\pm 1}, X_i^\pm$ 
subject to the relations
\begin{gather*}
K_i K_j = K_j K_i, \enspace
K_i X_j^\pm K_i^{-1} = 
q_i^{ \pm c_{ij}} X_j^\pm, \enspace
X_i^+ X_j^- - X_j^- X_i^+ = \delta_{ij} \frac{K_i - K_i^{-1}}{q_i - q_i^{-1}}, 
\\
\sum_{k=0}^{1-c_{ij}} (-1)^k
\begin{bmatrix}
1 - c_{ij} \\
k
\end{bmatrix}_{i}
(X_i^\pm)^{1 - c_{ij} -k} X_j^\pm (X_i^\pm)^k =0
\quad (\mbox{for} \, \,  i \neq j), 
\end{gather*}
where
\[
\quad [n]_{i} := \frac{q_i^n - q_i^{-n}}{q_i - q_i^{-1}}, \quad
[n]_i! := [1]_{i} \cdots [n]_{i}, \quad \mbox{and} \quad
\begin{bmatrix}
n \\
k
\end{bmatrix}_{i}
:= \frac{[n]_{i}}{ [k]_{i}! [n-k]_{i}!} \cdot
\]
It is a Hopf algebra with coproduct, antipode and counit given by 
\begin{gather*}
\Delta(K_i) = K_i \otimes K_i, \enspace 
\Delta(X_i^+) = X_i^+ \otimes 1 + K_i \otimes X_i^+, \enspace 
\Delta(X_i^-) = X_i^- \otimes K_i^{-1} + 1 \otimes X_i^-,
\\
S(K_i) = K_i^{-1}, \enspace  S(X_i^+) = - K_i^{-1} X_i^+,
S(X_i^-) = - X_i^- K_i, \enspace 
\ep(K_i) =1, \enspace \ep(X_i^\pm) = 0.
\end{gather*}
Let $U_q^\pm(\g)$ and $U_q^0(\g)$ be the 
subalgebras of $U_q(\g)$ generated by $\{X^\pm_i \mid 1 \leq i \leq r\}$ and $\{K_i^{\pm1} \mid 1 \leq i \leq r \}$, respectively. Denote
\[
U^{\geq}_q(\g) := U_q^0(\g) U_q^+(\g) \quad \mbox{and} \quad U^{\leq}_q(\g) := U_q^0(\g) U_q^-(\g).
\]

Let $\om$ be the involutive automorphism of $U_q(\g)$ given on generators by 
\begin{equation}
\label{om}
\om(X_i^\pm) = X_i^\mp, \enspace \om(K_i) = K_i^{-1}, \quad \forall 1 \leq i \leq r
\end{equation}
and $\tau$ be the involutive antiautomorphism of $U_q(\g)$ given on generators by
\begin{equation}
\label{tau}
\tau(X^\pm_i) = X^\pm_i,
\enspace
\tau(K_i) = K_i^{-1}, \quad
\forall\,  i \in [1,r].
\end{equation}
  
The algebra $U_q(\g)$ is $Q$-graded: 
\[
\deg K_i^{\pm 1} =0, \enspace \deg X_i^\pm = \pm \al_i, \quad \forall 1 \leq i \leq r.  
\]
For a graded subalgebra $A \subset U_q(\g)$, its graded components will be denoted by 
\[
A_\beta, \quad \beta \in Q.
\]
The degree (weight) of a homogeneous element $z\in A_\beta$
will be denoted by $\wt z:=\beta$.

The weight spaces of a finite-dimensional (left) $U_q(\g)$-module $V$ will be denoted by 
$$
V_\mu := \{ v \in V \mid K_i v = q^{ \lcor \mu, \al_i \rcor} v, \enspace 
\forall\,  i \in [1,r] \}, \quad \mu \in P.
$$
Such a module is called a \emph{type one} module if it is a direct sum of its weight spaces; those modules form a semisimple braided tensor category. The irreducible ones are parametrized by their highest weights which are dominant integral weights $P^+$, \cite[Theorem 5.10]{jantzen1996lectures}. Denote by $V(\mu)$ the irreducible module with highest weight $\mu \in P^+$. For $v \in V(\mu)$ and $\xi \in V(\mu)^*$, consider the matrix coefficient 
\begin{equation} 
\label{c-notation}
c_{\xi, v} \in U_q(\g)^* \quad 
\mbox{given by} \quad c_{\xi, v}(x) := \xi ( x v ),  \quad \forall x \in U_q(\g). 
\end{equation}
We fix a highest weight vector 
\[\
v_\mu \in V(\mu)_\mu.
\]

The quantized coordinate ring $R_q[G]$ of the connected simply connected algebraic group $G$ with Lie algebra $\g$ 
is the Hopf subalgebra of the restricted dual of $U_q(\g)$:
\[
R_q[G] = \bigoplus_{\mu \in P^+} \{ c_{\xi, v} \mid v \in V(\mu), \xi \in V(\mu)^* \}. 
\]
It is $P \times P$-graded with components 
\[
R_q[G]_{\nu, \nu'} = \{ c_{\xi, v} \mid
\mu \in \PP^+, \thickspace \xi \in (V(\mu)^*)_{\nu}, \thickspace v \in V(\mu)_{\nu'} \}, 
\quad \nu, \nu' \in P.
\]

Denote by $T_w$ (for $w \in W$) the Lusztig actions of the braid group $B_\g$ associated to $\g$ on $U_q(\g)$ \cite[\S 8.14]{jantzen1996lectures} and the finite-dimensional type one $U_q(\g)$-modules \cite[\S 8.6]{jantzen1996lectures}. Consider the extremal weight vectors, 
\[
v_{w \mu} := T^{-1}_{w^{-1}} v_\mu \in V(\mu)_{w\mu}, \quad 
\forall \mu \in \PP^+, w \in W.
\]
They only depend on $w \mu$, \cite[Proposition 39.3.7]{Lus:intro}. The extremal weight spaces $V(\mu)_{w\mu}$ are one dimensional. Let
\[
\xi_{w \mu} \in (V(\mu)^*)_{- w\mu}
\]
be the unique vector such that
$\lcor \xi_{w \mu}, T^{-1}_{w^{-1}} v_\mu \rcor =1$. 
The {\em{generalized quantum minors}} are given by
\begin{align}\label{eq:generalized-minor}
\De_{w \mu, u \mu} := c_{\xi_{w\mu}, v_{u\mu}} \in R_q[G]_{-w\mu, u \mu}, 
\quad \forall \mu \in P^+, w, u \in W.
\end{align}
They are equal to the Berenstein--Zelevinsky quantum minors \cite[Eq. (9.10)]{BerensteinZelevinsky05} defined by 
\[
\De_{w \mu, u \mu} (x) = \De_{\mu, \mu} ( T_{w^{-1}} \cdot x \cdot T_{u^{-1}}^{-1} ), \quad \forall\,  x \in U_q(\g).
\]
See \cite[\S 9.3]{GY13}
\subsection{Quantum double Bruhat cells}
\label{sec:d-Bruhat}
The subalgebras 
\begin{align*}
R^+ &= 
\Span \{c_{\xi, v} \mid \mu \in \PP^+, \thickspace
v \in V(\mu)_\mu, \thickspace \xi \in V(\mu)^* \} \subset R_q[G]
\quad \mbox{and}
\\
R^- & = \Span \{ c_{\xi, v} \mid \mu \in \PP^+, \thickspace
v \in V(\mu)_{w_0 \mu}, \thickspace \xi \in V(\mu)^* \} \subset R_q[G]
\end{align*}
generate $R_q[G]$, and more precisely, 
\[
R_q[G]= R^+ R^- = R^- R^+,
\]
see \cite[Proposition 9.2.2]{joseph2012quantum}. Consider the lowest weight vectors 
\[
v_{-\mu} := v_{w_0(- w_0 \mu)}= T^{-1}_{w_0} v_{- w_0 \mu} \in V(-w_0 \mu)_{-\mu}, \quad \forall \mu \in P^+
\]
and denote 
\[
c^+_\xi := c_{\xi, v_\mu } \in R^+ \quad \mbox{and} \quad 
c^-_{\xi'} := c_{\xi', v_{-\mu} } \in R^-,
\quad \forall \mu \in P^+, \xi \in V(\mu)^*, \xi' \in V(- w_0 \mu)^*
\]
and the quantum minors
\[
\De_{- w \mu, -\mu}:= \De_{w w_0 (-w_0 \mu), w_0 (-w_0\mu)} =
c_{\xi_{-w \mu},v_{-\mu}}.
\]
The Demazure modules
\[
V_w^+(\mu) = U_q^+(\g) V(\mu)_{w \mu} \subset V(\mu) \quad \mbox{and} \quad 
V_u^-(\mu) = U_q^-(\g) V(- w_0 \mu)_{- u \mu} \subset V(- w_0 \mu)
\]
give rise to the $P \times P$-homogeneous, completely prime ideals 
\begin{align*}
&I_w^\pm := \Span \{ c^\pm_\xi \mid \xi \in V_w^\pm(\mu)^\perp, \thickspace \mu \in \PP^+ \} 
\subset R^\pm, \\
&I_{w, u} := I_w^+ R^- + R^+ I_u^- \subset R_q[G],
\end{align*}
see \cite[Propositions 10.1.8 and 10.3.5]{joseph2012quantum}. 
\begin{definition} 
\label{def:normal}
An element $\Delta$ of a (noncommutative) ring $A$ is called {\em{normal}} if 
\[
\Delta A = A \Delta.
\] 
\end{definition}
The multiplicative subsets 
\[
E^\pm_w := \{ \De_{\pm w \mu, \pm \mu} \mid \mu \in \PP^+ \} \subset R^\pm, \quad
E_{w,u} := q^\Zset E_w^+ E_u^- \subset R_q[G],
\]
consist of nonzero normal elements of $R^\pm/I_w^\pm$ and $R_q[G]/I_{w,u}$, 
respectively, because
\begin{align}
\label{n1}
\De_{w \mu, \mu} a &= q^{\lcor w \mu, \nu \rcor - \lcor \mu, \nu' \rcor} 
a \De_{w \mu, \mu} \mod I_{w}^+ R^-,  
\\
\label{n2}
\De_{-w \mu, - \mu} a &= q^{ \lcor w \mu, \nu \rcor - \lcor \mu, \nu' \rcor} 
a \De_{-w \mu, - \mu} \mod R^+ I_w^-, 
\end{align}
$\forall \mu \in P^+$, $w \in W$,
$\nu, \nu' \in P$, $a \in R_q[G]_{-\nu, \nu'}$, see \cite[Eqs. (2.22) and (2.23)]{yakimov2014spectra}. 
The localization 
\begin{equation}
\label{qdBr}
R_q[G^{w,u}] := (R_q[G]/I_{w,u})[E_{w,u}^{-1}]
\end{equation}
is called the \emph{quantized coordinate ring of the double Bruhat cell 
$G^{w,u} := B^+ w B^+ \cap B^- u B^-$}. The quantized coordinate 
ring of the {\em{reduced double Bruhat cell}} $G^{w,u}/H$ is the subalgebra 
\[
R_q[G^{w,u}/H]
\]
of $R_q[G^{w,u}]$, 
which is the span of all elements of degrees in $P \times \{0\}$.
\subsection{Quantum Schubert cells}
\label{sec:q-Schubert}
For $w \in W$, the {\em{quantum Schubert cell}} algebras are the algebras
\[
U_q^\pm[w] := U_q^\pm(\g) \cap T_w(U_q^\mp(\g)).
\]
Given a reduced expression $w = s_{i_1} \ldots s_{i_l}$,
define the roots and root vectors 
\begin{equation}
\label{eq:beta}
\beta_k := w_{\leq k-1} (\al_{i_k}) \quad \mbox{and} \quad X^\pm_{\beta_k} := T_{w_{\leq k-1}} (X^\pm_{i_k}), \quad \forall \quad 1 \leq k \leq l. 
\end{equation}
The algebras $U_q^\pm[w]$ have PBW bases in the generators $X_{\be_1}^\pm, \ldots, X_{\be_l}^\pm$.
They were introduced by De Concini--Kac--Procesi \cite{DKP:solvable} and Lusztig \cite{Lus:intro} 
as the subalgebras of $U_q^\pm(\g)$ 
generated by the root vectors $X_{\be_1}^\pm, \ldots, X_{\be_l}^\pm$. Berenstein and Greenstein \cite{BG} conjectured that the two definitions lead to the same algebras and this was proved by Kimura \cite{kimura2016remarks} and Tanisaki \cite{tanisaki2017modules}. The opposite algebras of $U_q^\pm[w]$ were called {\em{quantum unipotent cells}} 
by Gei\ss, Leclerc and Schr\"oer \cite{GeissLeclercSchroeer11} (for the statement that the two algebras in question are anti-isomorphic see the display under Eq. (5.10) in \cite{GY21} and Eq. (6.3) therein).

The involutive automorphism $\omega$ of $U_q(\g)$, given by \eqref{om}, restricts to an isomorphism
$\om : U^+_q[w] \stackrel{\simeq}{\longrightarrow} U^-_q[w]$ and, on root vectors, it acts by 
\begin{equation}
\label{eq:om-prop}
\om (X_{\be_k}^+) = (-1)^{ \langle \be_k - \al_{i_k}, \rho\spcheck \rangle}  q^{- \langle \be_k - \al_{i_k}, \rho \rangle} X_{\be_k}^-, \quad \forall 1 \leq k \leq l.
\end{equation}

The algebras $R^\pm_w := (R^\pm/I^\pm_w)[(E^\pm_w)^{-1}]$ are $P \times P$-graded. Denote by $S_w^\pm$ their subalgebras consisting of elements of degrees in $P \times \{ 0 \}$; that is, $S_w^\pm$ are 
$Q \simeq Q \times \{ 0 \}$-graded. Their graded components will be denoted by $(S_w^\pm)_\beta$ for $\beta \in \QQ$; $(S_w^\pm)_0 = \C(q)$. The elements of $S^+_w$ and $S^-_u$ can be written in a nonunique way as
\[
c^+_\xi \De_{w \mu, \mu}^{-1} \quad \mbox{and} 
\quad \De_{- w \mu, - \mu}^{-1} c^-_{\xi'}
\]
for some $\mu \in P^+$, $\xi \in V(\mu)^*$ and 
$\xi' \in V(- w_0 \mu)^*$, cf. \cite[\S 10.3.1]{joseph2012quantum}.
(We omit coset notations for brevity.) 

The quantum quasi $R$-matrix $\RR^w$ is the element 
\begin{multline}
\label{Rw}
\RR^w := \sum_{m_1, \ldots, m_N \, \in \, \Zset_{\geq 0}}
\left( \prod_{j=1}^N
\frac{ (q_{i_j}^{-1} - q_{i_j})^{m_j}}
{q_{i_j}^{m_j (m_j-1)/2} [m_j]_{q_{i_j}}! } \right) 
\\
(X^+_{\be_N})^{m_N} \ldots (X^+_{\be_1})^{m_1} \otimes  
(X^-_{\be_N})^{m_N} \ldots (X^-_{\be_1})^{m_1}
\end{multline}
of the completion of $U^+_q[w] \otimes U^-_q[w]$ with respect to the descending filtration \cite[\S 4.1.1]{Lus:intro}. Denote the flip of the two components of $\RR^w$ by
$\RR^w_{\opp}$. Define the quasi $R$-matrices of $U_q(\g)$:
\[
\RR := \RR^{w_\ci}, \quad \RR_{\opp} := \RR^{w_\ci}_{\opp}. 
\]
For a linear map $\theta : U_q(\g) \to U_q(\g)$, define
\[
{}^{\theta \times 1}\RR:= (\theta \otimes \id)\RR,
\]
and similarly,
${}^{\theta \times 1}\RR^w$, ${}^{\theta \times 1}\RR_{\opp}$,
${}^{1 \times \theta}\RR$.

\bth{isom} \cite[Theorem 2.6]{yakimov2014spectra}
For all complex simple Lie algebras $\g$ and $w \in W$, 
the maps $\vp^\pm_w : S^\pm_w \to U^\mp_q[w]$ given by
\begin{align*}
& \vp^+_w \big( c^+_\xi \De_{w\mu, \mu}^{-1} \big) 
:= \big( c_{\xi, v_{w \mu}} \otimes \id \big)
{}^{\tau \times 1} \RR^w, \enspace 
\forall\,  \mu \in \PP^+, \thickspace \xi \in V(\mu)^*,
\\
& \vp^-_w \big( \De_{-w\mu, -\mu}^{-1} c^-_\xi  \big) 
:= \big( \id \otimes c_{\xi, v_{- w \mu}} \big)
{}^{1 \times \tau} \RR^w_{\opp}, \enspace 
\forall\,  \mu \in \PP^+, \thickspace \xi \in V(- w_0\mu)^*
\end{align*}
are well defined $Q$-graded algebra anti-isomorphisms and isomorphisms, 
respectively {\em{(}}recall \eqref{tau}{\em{)}}. 
In both formulas, one can replace $\RR^w$ with $\RR$ leading to the same maps $\vp^\pm_w$. 
\eth

The quantum double Bruhat cell algebras $R_q[G^{w,u}]$ can be reconstructed from the algebras $S_w^+$ and $S_u^-$, and in light of Theorem \ref{tisom}, from the algebras $U_q^-[w]$ and $U_q^+[u]$ as well. Joseph proved \cite[\S 10.3.2]{joseph2012quantum} that there exists a (unique) $Q$-graded algebra $S^+_w \bowtie S^-_u$ such that
\begin{enumerate}
\item[(i)] $S^+_w \bowtie S^-_u \simeq S^+_w \otimes_{\KK} S^-_u \simeq S^-_u \otimes_{\KK} S^+_w $ as $\KK$-vector spaces and
the canonical embeddings of $S^+_w$ and $S^-_u$ in it are graded algebra embeddings, and
\item[(ii)] the elements of $S^+_w$ and $S^-_u$ satisfy the commutation relations
\begin{align}
\label{commRR}
&\big[ \De_{- u \mu', -\mu'}^{-1} c^-_{\xi'}  \big] 
\big[ c^+_\xi \De_{w \mu, \mu}^{-1} \big]
= q^{ - \lcor \nu' + u \mu' , \nu - w \mu  \rcor} 
\big[ c^+_\xi \De_{w \mu, \mu}^{-1} \big]
\big[ \De_{ - u \mu', - \mu'}^{-1} c^-_{\xi'}  \big]
\\
\nn
+ & \sum_{i}
q^{ - \lcor \nu' + \ga_i + u \mu', \nu - \ga_i - w \mu \rcor } 
\big[ c^+_{S^{-1}(x_i) \xi} \De_{w \mu, \mu}^{-1} \big] .
\big[ \De_{-u \mu', - \mu'}^{-1} c^-_{S^{-1}(x'_{i}) \xi'} \big],
\end{align}
for all $\mu, \mu' \in \PP^+$, $\xi \in( V(\mu)^*)_{-\nu}$, $\xi' \in (V(- w_0 \mu')^*)_{- \nu'}$, where $\RR = \sum_i x_i \otimes x'_i$, and $x_i^\pm \in U_q^\pm(\g)_{\pm \ga_i}$. 
\end{enumerate}
\begin{eg} 
\label{sl3-1}
Consider the case of $G = SL_3$ and $w = u = w_0 \in S_3$. In this case 
\[
S^+_{w_0} \bowtie S^-_{w_0} \simeq 
U_q^-[w_0]^{\opp} \bowtie U_q^+[w_0] \simeq 
U_q^-(\mathfrak{sl}_3)^{\opp} \bowtie U_q^+(\mathfrak{sl}_3). 
\]
This algebra is isomorphic to Kashiwara's \cite[Sect. 3.3]{Kas:crystal} bosonic algebra $\mathscr{B}_q(\mathfrak{sl}_3)$, see \cite[Remark 4.8]{goodearl2016berenstein}. It is presented as the algebra with generators $X_1^\pm$, $X_2^\pm$ and relations
\begin{gather*}
(X_i^\pm)^2 X_{3-i}^{\pm} - (q + q^{-1}) X_i^\pm X_{3-i}^\pm X_i^\pm + X_{3-i}^\pm (X_i^\pm)^2=0,  
\\
X_i^+ X_{3-i}^- = q^{-1} X_{3-i}^- X_i^+,  \\ 
X_i^+ X_i^- = q^2 X_i^- X_i^+ + (q^{-1} - q)^{-1}. 
\end{gather*}
The algebra has an action of $(\KK^\times)^2$ by algebra automorphisms given by
\[
(s_1, s_2) \cdot X_i^\pm := s_i^{\pm 1} X_i^\pm, \quad i =1,2 
\]
for $s_1, s_2 \in \KK^\times$. Consider the elements 
\[
X_{12}^+:= X_1^+ X_2^+ - q^{-1} X_2^+ X_1^+ \quad \mbox{and} \quad X_{12}^-:= X_1^- X_2^- - q X_2^- X_1^-.
\]
 The 
six elements 
\[
X_1^\pm,\ X_2^\pm,\ X_{12}^\pm
\]
satisfy the relations
\[
\begin{gathered}
\begin{aligned}
X_2^+ X_1^+ &= q X_1^+ X_2^+ - q X_{12}^+, \qquad  &X_1^- X_2^- &= q X_2^- X_1^- + X_{12}^-, \\
X_{12}^+ X_1^+ &= q^{-1} X_1^+ X_{12}^+,  &X_2^+ X_{12}^+ &=  q^{-1} X_{12}^+ X_2^+, \\
X_{12}^- X_2^- &= q^{-1} X_2^- X_{12}^-,  &X_1^- X_{12}^- &=  q^{-1} X_{12}^- X_1^-, \\
X_{12}^+ X_1^- &= q X_1^- X_{12}^+,  &X_{12}^+ X_2^- &= q X_2^- X_{12}^+ - q^{-1} X_1^+, \\
X_1^+ X_{12}^-  &= q X_{12}^- X_1^+,  &X_2^+ X_{12}^-  &= q X_{12}^- X_2^+  + X_1^-,
\end{aligned}  \\
X_{12}^+ X_{12}^- = q^2 X_{12}^- X_{12}^+ + (q^2 -1) X_1^- X_1^+ + (q^{-1} - q)^{-1},
\end{gathered}
\]
see \cite[Example 2.4]{goodearl2016berenstein}. These relations, together with the $(\KK^\times)$-action turn the algebra into a symmetric CGL extension \cite[Definitions 3.3 and 3.12]{GY13}
in the order of adjoining the generators $X_2^-,\ X_{12}^-,\ X_1^-,\ X_1^+,\ X_{12}^+,\ X_2^+$.
\end{eg}
The compositions 
\begin{align*}
&S^+_w \hra (R^+/I_w^+)[(E^+_w)^{-1}] \hra (R_q[G]/I_{w,u})[E_{w,u}^{-1}]=R_q[G^{w,u}] 
\quad \mbox{and} \quad \\
&S^-_u \hra (R^-/I_u^-)[(E^-_u)^{-1}] \hra (R_q[G]/I_{w,u})[E_{w,u}^{-1}]=R_q[G^{w,u}].
\end{align*}
give rise to an embedding 
\begin{equation}
\label{SRembed}
(S^+_w \bowtie S^-_u) [p_i^{-1}, i \in \cS(w,u)] \hookrightarrow  R_q[G^{w,u}/H].
\end{equation}

Fix an isomorphism of $U_q(\g)$-modules $V(\vpi_i)^* \simeq V(- w_0 \vpi_i)$ such that $\xi_{\vpi_i} \mt v_{- \vpi_i}$. 
For a pair $\{\xi_{i, j} \}$ and $\{\xi'_{i,j} \}$
of dual bases of $V(\vpi_i)^* \simeq V(- w_0 \vpi_i)$ and $V(- w_0 \vpi_i)^*$, consider the elements 
\[
p_i:= \sum_j \big[ \De_{- u \vpi_i, - \vpi_i}^{-1} c^-_{\xi'_{i,j}}  \big]
\big[ c^+_{\xi_{i,j}} \De_{w \vpi_i, \vpi_i}^{-1} \big]  \in (S^+_w \bowtie S^-_u)_{(w-u)\vpi_i}
\]
(this is independent of the choice of dual bases). These elements are nonzero normal elements \cite[Theorem 5.1]{goodearl2016berenstein} and are 
non-scalars if and only if $i \in \cS(w,u)$, \cite[Proposition 7.6]{yakimov2014spectra}.

In the embedding \eqref{SRembed},
\[
p_i \mt \De_{- u \vpi_i, - \vpi_i}^{-1} \De_{w \vpi_i, \vpi_i}^{-1},
\]
\cite[Eq. (3.19)]{goodearl2016berenstein} and \eqref{SRembed} gives rise to the following realization of the quantized coordinate ring of the reduced double Bruhat cell $G^{w,u}/H$:
\begin{equation}
\label{eq:GwuSwu}
(S^+_w \bowtie S^-_u) [p_i^{-1}, i \in \cS(w,u)]
\simeq  R_q[G^{w,u}/H],
\end{equation}
see \cite[Eq. (3.23)]{goodearl2016berenstein}.

Let $\LL_w^+$ and $\LL_u^-$ be the subalgebras of $R_q[G^{w,u}]$ generated by the elements of $E_w^+ \cup (E_w^+)^{-1}$ 
and $E_u^- \cup (E_u^-)^{-1}$. They are isomorphic to Laurent polynomial rings in the $r$ indeterminates
\[
\De_{w \vpi_i, \vpi_i}^{ \pm 1} \quad \mbox{and} \quad
\De_{-u \vpi_i, - \vpi_i}^{ \pm 1}, \quad 1 \leq i \leq r. 
\]
By \cite[10.3.2(5)]{joseph2012quantum}, 
\begin{align}
\label{dBisom}
R_q[G^{w,u}] &\simeq
\big( (S^+_w \bowtie S^-_u)[p_i^{-1}, i \in \cS(w,u)] \big) \# \LL_w^+
\\
&\simeq
\big( (S^+_w \bowtie S^-_u)[p_i^{-1}, i \in \cS(w,u)] \big) \# \LL_u^-,
\nn
\end{align}
where the smash products are defined by using the commutation relations \eqref{n1}-\eqref{n2}. 

\subsection{Dual canonical basis}

Choose and fix a reduced expression $s_{i_1}\cdots s_{i_{l(w)}}$ for $w$. For any $c\in \N$, let us denote $\PBW(c \beta_k):=\frac{1}{[c]_{i_k}!}(X^{-}_{\beta_k})^{c}$. Define the PBW basis elements
\begin{align*}
    \PBW(\uc):=\PBW(c_l \beta_{l(w)})\cdots \PBW(c_1 \beta_1), \forall \uc\in \N^{[1,l(w)]}.
\end{align*}

Let $(\ ,\ )_L$ denote Lusztig's pairing on $U^-_q[w]$, which is a symmetric and non-degenerate $\K$-bilinear form. By \cite[Proposition 4.19]{Kimura10}\cite[Proposition 38.2.1]{Lus:intro}, we have
\begin{align*}
    (\PBW(\uc),\PBW(\uc'))_L=\prod_{k=1}^{l(w)}\delta_{c_k,c'_k}\prod_{s=1}^{c_k}\frac{1}{1-q_{i_k}^{2s}}.
\end{align*}
Define the corresponding dual PBW basis elements.
\begin{align*}
\dPBW(\uc):=\dPBW(\sum_k c_k\beta_k):=&\frac{1}{( \PBW(\uc),\PBW(\uc))_L}\PBW(\uc).
\end{align*}
In particular, we obtain the dual PBW generators.
\begin{align*}
    \dPBW(\beta_k):=\dPBW_{\beta_k}:=(1-q_{i_k}^2)X^{-}_{\beta_k}.
\end{align*}
By \cite[Proposition 4.23]{Kimura10}, we have 
\begin{align}\label{eq:dual-PBW-multiple}
\dPBW(c \beta_k)\dPBW(c' \beta_k)=q^{\Hf \langle c \beta_k,c'\beta_k\rangle}\dPBW((c+c') \beta_k),\ \forall c,c'\in\N.
\end{align}

Define the dual bar involution $\sigma$ on $U^-_q[w]$ by $(\sigma(x),y)_L=\overline{(x,\overline{y})_L}$. Then we have $\sigma(q)=q^{-1}$, $\sigma(X^{-}_{\beta_k})=-q_{i_k}^2 X^{-}_{\beta_k}$, and $\sigma(\dPBW_{\beta_k})=\dPBW_{\beta_k}$. For homogeneous $x,y\in U^-_q[w]$, we have $\sigma(xy)=q^{\langle \wt x,\wt y\rangle}\sigma(y)\sigma(x)$, see \cite[Proposition 3.10]{kimura2017twist}.

The dual canonical basis element $\dCan(\uc)$ is defined as the unique $\sigma$-invariant element such that 
\begin{align}\label{eq:PBW-canonical}
\dCan(\uc)=\dPBW(\uc)+ \sum_{\uc'<_{\lex}\uc} b_{\uc'}\dPBW(\uc'),\quad b_{\uc'}\in q\Z[q].
\end{align}
They form the dual canonical basis $\dCan$ for $U^-_q[w]$, which was denoted as $B^{\up}_{-1}$ in \cite[Section 7.3]{qin2020dual}, \cite{Kimura10}. Note that \eqref{eq:PBW-canonical} is homogeneous in the $Q$-grading. The following result is well-known and can be checked directly.

\begin{lem}\label{lem:dCan-inv}
The dual canonical basis $\dCan$ for $U^-_q[w]$ is independent of the choice of the reduced word for $w$.
\end{lem}

Finally, choose and fix any $\Z$-valued linear function $\psi$ on $Q$. As in \cite[Section 7.1]{qin2023analogs} for any homogeneous $x\in U^-_q[w]$, define \begin{align*}
    \sigma'_{\psi}(x)=q^{-\Hf\langle \wt x,\wt x\rangle +\psi(\wt x)}\sigma(x)
\end{align*}
and extend it to $U^-_q[w]$ by additivity. The map $\sigma'_{\psi}$ is a skew linear automorphism $U^-_q[w]$, i.e., 
$\sigma'_{\psi}(xy)=\sigma'_{\psi}(y)\sigma'_{\psi}(x)$, $\forall x,y\in U^-_q[w]$.

Further, we consider the $\Q(q^{\Hf})$-linear automorphism $\cor_{\psi}$ on $U^-_q[w]$ such that, for any homogeneous $x\in U^-_q[w]_\gamma$, 
\begin{align}\label{eq:cor}
    \cor_{\psi} x:= q^{-\frac{1}{4}\langle \wt x,\wt x\rangle +\Hf \psi(\wt x)}x.
\end{align}
Here, $\cor$ abbreviates $q$-factor correction. Then $\sigma(x)=x$ if and only if $\sigma'_{\psi}(\cor_{\psi} x)=\cor_{\psi} x$. We will abbreviate $\sigma'_{\psi}$ by $\sigma'$ and $\cor_{\psi}$ by $\cor$ for simplicity from now on.

Note that, for homogeneous $x,y$, we have
\begin{align}\label{eq:cor-prod}
    \cor (xy)=q^{-\Hf\langle \wt x,\wt y\rangle}\cor(x)\cor(y).
\end{align}
Combining with \eqref{eq:dual-PBW-multiple}, we deduce the following.
\begin{align}\label{eq:cor-dual-PBW}
    \begin{split}
        \cor(\dPBW(c \beta_k))&=(\cor \dPBW_{\beta_k})^c,\ \forall c\in\N.\\
        \cor(\dPBW(\uc))&=q^{-\Hf\sum_{k>j}\langle c_k \beta_k, c_j\beta_j\rangle}(\cor \dPBW_{\beta_{l(w)}})^{c_{l(w)}} \cdots (\cor \dPBW_{\beta_1})^{c_1},\ \forall \uc\in\N^{[1,l]}.  
    \end{split}
\end{align}

\section{Cluster algebra structures on quantum double Bruhat cells}
\label{sec:dBruhat}
In this section we prove Theorem A. We also obtain several results on common triangular bases of (reduced) quantum double Bruhat cells that are of independent interest. The quantum groups and quantum cluster algebras that are considered are defined over $\K:=\C(q^\Hf)$. 

\subsection{Quantum seeds on reduced double Bruhat cells}

Fix two Weyl group elements $w,u$. Choose a reduced word $\uzeta$ for $w$ and a reduced word $\ueta$ for $u$, whose letters belong to $[1,r]$. Choose a signed word $\ubi$ such that it can be obtained from the signed word $(-\uzeta,\ueta)$ by left reflections and flips. Denote $\rsd=\rsd(\ubi)$, $\dsd=\dsd(\ubi)$. Then $l:=l(\ubi)$ equals $l(w)+l(u)$. Recall that $\dI=[1,l]$ and $\ddI=\{-i \mid i\in [1,r]\}\sqcup [1,l]$. We introduce the following elements depending on the choice of $\ueta$ and $\uzeta$:
\begin{align}\label{eq:subword}
    (w_{>k})^{-1}&=s_{\zeta_{l(w)}}\cdots s_{\zeta_{k+1}},\quad \forall k\in [1,l(w)],\\
    u_{\leq k}&=s_{\eta_1}\cdots s_{\eta_k},\quad \forall k\in [1,l(u)].
\end{align}
In \cite[Eq. (7.2)]{goodearl2016berenstein}, $(w_{>k})^{-1}$ was denoted by $w^{-1}_{<k}$. 

Denote the distance function
\begin{align*}
    d(k):= | k - l(w) - 1/2| + 1/2 =
    \begin{cases}
    l(w)+1-k, &k\in [1,l(w)]\\
    k-l(w),& k\in [l(w)+1,l].
    \end{cases}
\end{align*}

Recall that $\beta_k$ denotes the $k$-th root associated with the reduced word $\uzeta$ of $w$, see
\eqref{eq:beta}. Similarly, let $\beta'_k$ denote the $k$-th root associated with the reduced word $\ueta$ of $u$, and $X_{\beta'_k}^{\pm}$ the associated root vectors. Recall that we have the anti-isomorphism $\varphi^+_w:S^w_+\simeq U^-_q[w]$ and the isomorphism $\varphi^-_u:S^u_-\simeq U^+_q[u]$, see Theorem \ref{tisom}. Denote $C_{\beta_k}^+:=(\varphi^+_w)^{-1}(X_{\beta_k}^-)$ and $C_{\beta'_k}^-:=(\varphi^-_u)^{-1}(X_{\beta'_k}^+)$.

An initial quantum seed $\sd^{\GY}$ on $S^+_w \bowtie S^-_u$ was introduced in \cite[Theorem 6.2(c)]{goodearl2016berenstein}. It can be mutated to another quantum seed $\sd^{\GY}_{\id}$, where we choose the permutation $\sigma$ to be the identity in \cite[Theorem 6.2(a)]{goodearl2016berenstein}. We have $\tB(\sd^{\GY})=-\tB(\rsd(-\uzeta,\ueta))$ and $\tB(\sd^{\GY}_{\id})=-\tB(\rsd(\uzeta^\op,\ueta))$, see Example \ref{sl3-2}.

The $\Lambda$-matrix of $\sd^{\GY}$ will be recalled later. We 
set 
\begin{equation}
\label{eq:minus}
\Lambda(\rsd):=-\Lambda(\sd^{\GY}).
\end{equation}
Because of the minus sign, we obtain the $\K$-algebra anti-isomorphism
\begin{align*}
    \rkappa:\bUpClAlg(\rsd)\simeq S^+_w \bowtie S^-_u.
\end{align*}
Note that $\sd^{\GY}=\rsd^\op$, where $\rsd^\op$ is the opposite seed of $\rsd$.

The anti-isomorphism $\rkappa$ establishes a cluster structure on the algebra $S^+_w \bowtie S^-_u$. Let us describe $\rkappa$ explicitly. The elements $\rkappa x_k(\rsd)$ are denoted by $\overline{y}_{[k^{\min},k]}$ in \cite[Theorem 6.1, Proposition 9.4]{goodearl2016berenstein}.  Recall that the mutation sequence $\Sigma$ is a green-to-red sequence. In this case, $\Sigma$ starts with the initial seed $\sd^{\GY}$ and ends with the seed of 
$S^+_w \bowtie S^-_u$ associated to the CGL extension presentation of this algebra for the reverse order of the PBW generators \cite[Theorem 8.2(b)]{GY13} (this is the seed corresponding to $\tau$ equal to the longest element of the symmetric group $S_{l(w) + l(u)}$). The sequence is constructed by successively pulling to the right one of the PBW generators as in the diagram after Lemma 5.5 in  \cite{GY13} and considering the corresponding seeds as in \cite[Theorem 8.2(b)]{GY13}; the mutations are provided by 
\cite[Theorem 8.2(c)]{GY13}. We next describe 
the explicit form of all interval variables $W_{[j,k]}$, where $|\bi_{j}|=|\bi_{k}|$ equals some $i\in [1,r]$. Their images $\rkappa W_{[j,k]}$ are denoted by $\overline{y}_{[j,k]}$ in \cite[Theorem 6.1]{goodearl2016berenstein}\cite[Theorem 10.1]{GY13}, 
and are given as follows:
\begin{itemize}
    \item For $1\leq j\leq k\leq l(w)$, \begin{align*}
        \rkappa W_{[j,k]}=\sqrt{q}^{\Hf|| (w_{\leq d(j)}-w_{<d(k)})\varpi_i ||^2}(\varphi^+_w)^{-1}\big((\Delta_{w_{<d(k)}\varpi_i,w_{\leq d(j)}\varpi_i}\otimes \id)^{\tau\otimes 1}\RR^w\big).
    \end{align*}
    \item For $l(w)< j\leq k\leq l$, \begin{align*}
        \rkappa W_{[j,k]}=\sqrt{q}^{\Hf|| (u_{\leq d(k)}-u_{<d(j)})\varpi_i ||^2}(\varphi^-_u)^{-1}\big((\Delta_{\varpi_i,\varpi_i}\otimes \id)\big((T_{u_{\leq d(k)}^{-1}}\otimes \id)\cdot ^{S\tau\otimes 1}\RR^u_{\op} \cdot (T_{u_{< d(j)}}\otimes \id) \big)\big).
    \end{align*}
    \item For $1\leq j\leq l(w)< k\leq l$, \begin{align*}
    \rkappa W_{[j,k]}=\sqrt{q}^{\Hf|| (u_{\leq d(k)}-w_{<d(j)})\varpi_i ||^2}(\varphi_{w,u})^{-1}\big((\Delta_{\varpi_i,w_{\leq d(j)}\varpi_i}\otimes \id)\big((T_{u_{\leq d(k)}^{-1}}\otimes \id)\cdot ^{S\tau\otimes 1}\RR^u_{\op} {^{\tau\otimes 1}\RR^w} \big)\big).
    \end{align*}
\end{itemize}
Here, $w_{\leq d(j)},w_{<d(k)},u_{\leq d(k)}$ are defined as in \eqref{eq:subword} and we denote $\varphi_{w,u}:=\varphi_{w}\otimes \varphi_{u}$ as in \cite[(4.18)]{goodearl2016berenstein}. We note that \cite[Theorem 6.1]{goodearl2016berenstein} states this result for $j<k$, but the same formula holds for $j=k$ as in \cite[Theorem 10.1]{GY13}. In particular, the image of $W_k=W_{[k,k]}$ satisfies the following equations, see \cite[Eq. (6.2)]{goodearl2016berenstein}:
\begin{align*}
    \rkappa W_k&=q_{i}^{\Hf}(q_{i}^{-1}-q_{i})C^+_{\beta_{d(k)}}, \quad \forall k\in[1,l(w)],\\
    \rkappa W_{k}&=q_{i}^{-\Hf}(q_{i}^{-1}-q_{i})q^{2\langle \beta'_{d(k)}, \rho\rangle}C^-_{\beta'_{d(k)}},\quad \forall k\in[1+l(w),l].
\end{align*}
Applying $\varphi^+_w$ and $\varphi^-_u$ to these equalities, respectively, we have 
\begin{align*}
    \varphi^+_w \rkappa W_k&=q_{i}^{\Hf}(q_{i}^{-1}-q_{i})X^-_{\beta_{d(k)}}=q^{-\frac{1}{4}\langle \beta_{d(k)},\beta_{d(k)}\rangle}\dPBW_{\beta_{d(k)}},\quad \forall k\in[1,l(w)],\\
    \varphi^-_u\rkappa W_{k}&=q_{i}^{-\Hf}(q_{i}^{-1}-q_{i})q^{2\langle \beta'_{d(k)}, \rho\rangle}X^+_{\beta'_{d(k)}},\quad \forall k\in[1+l(w),l].
\end{align*}

Following \cite[(6.9)]{goodearl2016berenstein}, we introduce a skew-symmetric bilinear form $\nu$ on $\Z^{[1,l]}$ whose matrix in the standard basis $\{ f_1, \ldots, f_l \}$ has entries
\begin{align*}
    \nu_{kj}&=-\langle \beta_{d(k)}, \beta_{d(j)} \rangle,\quad j<k\in [1,l(w)],\\
    \nu_{kj}&=-\langle \beta'_{d(k)}, \beta'_{d(j)} \rangle,\quad j<k\in [l(w)+1,l],\\
    \nu_{kj}&=\langle \beta'_{d(k)}, \beta_{d(j)} \rangle,\quad j\leq l(w)<k.
\end{align*}
Denote $f_{[k^{\min},k]}:=f_{k^{\min}}+f_{k^{\min}[1]}+\cdots+f_{k}$ for $k\in \dI$. Then entries of the $\Lambda$-matrix of the seed 
$\sd^{\GY}$ are given by
\begin{align*}
    \Lambda(\sd^{\GY})_{kj}=\nu(f_{[k^{\min},k]},f_{[j^{\min},j]}), &\quad j,k\in [1,l].
\end{align*}

Recall that we have the isomorphism $R_q[G^{w,u}/H]\simeq (S^+_w \bowtie S^-_u) [p_i^{-1}, i \in \cS(w,u)]$, see \eqref{eq:GwuSwu}. Identifying these two algebras and taking the localization at the frozen variables, we extend $\rkappa$ to the following anti-isomorphism, see \cite[Theorem 6.2(d)]{goodearl2016berenstein}:
\begin{align*}
    \rkappa:\upClAlg(\rsd)\simeq R_q[G^{w,u}/H].
\end{align*}
If $z\in R_q[G^{w,u}/H]$ is homogeneous of weight $(\gamma,0)$, we say that the weight of $z$ is $\wt z:=\gamma$. We will view $\upClAlg(\rsd)$ as a $P$-graded algebra by defining the weight of $\rkappa^{-1}z$ to be $\wt z$. We have, 
\begin{align*}
    \wt W_k&=-\beta_{d(k)},\quad \forall k\in[1,l(w)],\\
    \wt W_{k}&=\beta'_{d(k)},\quad \forall k\in[1+l(w),l],\\
    \wt W_{[j,k]}&=\wt W_j+\wt W_{j[1]}+\cdots +\wt W_{k},\quad \forall j,k\in[1,l], 
    |\bi_{j}|=|\bi_{k}|.
\end{align*}
It follows that
\begin{align}\label{eq:Lambda-GY}
    \begin{split}
        \nu_{kj}&=-\langle \wt W_k,\wt W_j\rangle,\ \forall k>j,\\
        \Lambda(\rsd)_{kj}&=-\Lambda(\sd^{\GY})_{kj}=-\nu(\wt W_{[k^{\min},k]},\wt W_{[j^{\min},j]}),\ \forall j,k\in [1,l],\\
        \Lambda(\rsd)(\deg W_k,\deg W_j)&=-\Lambda(\sd^{GY})(\deg W_k,\deg W_j)=-\nu_{kj},\  \forall j,k\in [1,l].      
    \end{split}
\end{align}

\begin{eg}
\label{sl3-2} We continue Example \ref{sl3-1}, and
consider the case of $G = SL_3$, $w=u =w_0 \in S_3$. Choose the reduced words $\uzeta=\ueta=(1,2,1)$ for $w$ and $u$, respectively.

The initial quantum seed $\sd^{\GY}_{\id}$ from \cite{goodearl2016berenstein} for $S^+_{w_0} \bowtie S^-_{w_0}$ , the algebra in Example \ref{sl3-1}, has cluster variables $x_k$, $k\in[1,6]$:
\[
\begin{aligned}
x_1=\rkappa W_{[1,1]} &= q^{\Hf} t^{-1} X_2^-, \qquad x_2=\rkappa W_{[2,2]} = q^{\Hf} t^{-1} X_{12}^-,\\
x_3=\rkappa W_{[1,3]}&= q^{\frac{3}{2}} t^{-2} (X_2^- X_1^- - t X_{12}^-), 
\\
x_4=\rkappa W_{[1,4]} &= q^{\frac{5}{2}} t^{-3} (X_2^- X_1^- X_1^+ - t X_{12}^- X_1^+ - q^{-1} t^2 X_2^-), \\
x_5=\rkappa W_{[2,5]}&= - q^3 t^{-2} (X_{12}^- X_{12}^+ + X_1^- X_1^+ - q^{-1} t^2), \\
x_6=\rkappa W_{[1,6]} &= q^4t^{-4} (X_2^- X_1^- X_1^+ X_2^+ - t X_{12}^-X_1^+ X_2^+ - q^{-1} t^2 X_2^- X_2^+  \\
&\hphantom{= X_2^- X_1^- X_1^+ X_2^+} \ + q t X_2^- X_1^- X_{12}^+ - q t^2 X_{12}^- X_{12}^+ + q^{-2} t^4),
\end{aligned}
\] 
where
\[
t: = (q^{-1} - q)^{-1}.
\]
The frozen cluster variables are $\rkappa W_{[2,5]}$ and $\rkappa W_{[1,6]}$. Note that we have $d(k)=4-k$ for $k\in [1,3]$. We can compute the matrices
\begin{align*}
(\nu_{jk})=\begin{pmatrix}
0&1&-1&1&-1&-2\\
-1&0&1&-1&-2&-1\\
1&-1&0&-2&-1&1\\
-1&1&2&0&1&-1\\
1&2&1&-1&0&1\\
2&1&-1&1&-1&0
\end{pmatrix},\ \Lambda=\begin{pmatrix}
0&1&-1&0&0&-2\\
-1&0&0&-1&-2&-2\\
1&0&0&-1&-2&-2\\
0&1&1&0&0&-2\\
0&2&2&0&0&0\\
2&2&2&2&0&0
\end{pmatrix}.
\end{align*}
The exchange matrix of the seed is the adjacency matrix of the ice quiver 
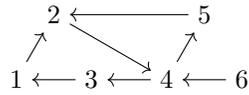
\[
\begin{tikzpicture}[
  vc/.style={circle,    draw, inner sep=1.5pt, minimum size=4mm},
  vs/.style={rectangle, draw, inner sep=1.5pt, minimum size=4mm}
]
\node[vc] (1) at (0,0) {1};
\node[vc] (3) at (1,0) {3};
\node[vc] (4) at (2,0) {4};
\node[vs] (6) at (3,0) {6};
\node[vc] (2) at (0.5,1.732/2) {2};
\node[vs] (5) at (2.5,1.732/2) {5};
\draw[<-] (1) -- (3);
\draw[<-] (3) -- (4);
\draw[<-] (4) -- (6);
\draw[<-] (2) -- (5);
\draw[<-] (2) -- (1);
\draw[<-] (4) -- (2);
\draw[<-] (5) -- (4);
\end{tikzpicture}
\]
where mutable vertices are circled and frozen boxed.
This exchange matrix equals $-\tB(\rsd(\uiota))$, where $\uiota=(\uzeta^\op,\ueta)=(1,2,1,1,2,1)$. Therefore, we can identify the opposite seed $\rsd(\uiota)^\op$ with $\sd^{\GY}_{\id}$.

Two nontrivial interval variables $\rkappa W_{[j,k]}$, $j\neq k$, that do not already appear in the initial seed are 
\begin{gather*}
\rkappa W_{[3,4]} = q t^{-2} (X_1^- X_1^+ - q^{-1} t^2),
\\
\rkappa W_{[3,6]} = q^{\frac{3}{2}} t^{-3}(X_1^- X_1^+ X_2^+ - q^{-1} t^2 X_2^+ + q t X_1^- X_{12}^+). 
\end{gather*}

Note that $\uiota$ and $(-\uzeta,\ueta)$ are related by left reflections and a flip. Indeed, the seed $\mu_1 \rsd(\uiota)^\op$ equals $\rsd(-\uzeta,\ueta)^\op$. The exchange matrix of the seed $\rsd(-\uzeta,\ueta)^\op$ corresponds to the following ice quiver
\[
\begin{tikzpicture}[
  vc/.style={circle,    draw, inner sep=1.5pt, minimum size=4mm},
  vs/.style={rectangle, draw, inner sep=1.5pt, minimum size=4mm}
]
\node[vc] (1) at (0,0) {1};
\node[vc] (3) at (1,0) {3};
\node[vc] (4) at (2,0) {4};
\node[vs] (6) at (3,0) {6};
\node[vc] (2) at (0.5,1.732/2) {2};
\node[vs] (5) at (2.5,1.732/2) {5};
\draw[->] (1) -- (3);
\draw[->] (3) -- (2);
\draw[<-] (3) -- (4);
\draw[<-] (4) -- (6);
\draw[<-] (2) -- (5);
\draw[->] (2) -- (1);
\draw[<-] (4) -- (2);
\draw[<-] (5) -- (4);
\end{tikzpicture}
\]
It coincides with the exchange matrix of $\sd^{\GY}$, see \cite[Example 2.16]{goodearl2016berenstein}. By \cite{goodearl2016berenstein}, the cluster variable $x_1':=x_1(\rsd(-\uzeta,\ueta)^\op)$ equals $\rkappa W_{[3,3]}=q^{\frac{1}{2}}t^{-1}X_1^-$. It is straightforward to check that $x_1 x_1'=q^{\frac{1}{2}} x_2 + q^{-\frac{1}{2}} x_3 $.
\end{eg}

\subsection{Quantum seeds on double Bruhat cells}
Take $\ubi=(-\uzeta,\ueta)$. 
Denote $\rsd=\rsd(\ubi)$ and $\dsd=\dsd(\ubi)$. Recall that $R_q[G^{w,u}]$ is $P\times P$-graded and that $R_q[G^{w,u}/H]$ is spanned by its graded components of degrees in $P\times \{0\}$. In \cite{goodearl2016berenstein}, an anti-isomorphism
\begin{align*}
    \dkappa:\upClAlg(\dsd)\simeq R_q[G^{w,u}].
\end{align*}
was defined, whose definition will be postponed to \eqref{eq:dkappa-def}. The ``anti'' part is because of the minus sign in \eqref{eq:minus}; note the interchange of the order of the two Weyl group elements. 

Via the map $\dkappa$, $\upClAlg(\dsd)$ becomes a $P\times P$-graded algebra. By \cite{goodearl2016berenstein}, we can identify $\upClAlg(\rsd)$ with the subalgebra spanned by the $P\times \{0\}$-graded components of $\upClAlg(\dsd)$, such that $x_k(\rsd)$ is identified with $\factor_k\cdot x_k(\dsd)$, where $\factor_k$ is a Laurent monomial of $x_{-i}(\dsd)$, $i\in[1,r]$.  Under this identification, we have
\begin{align*}
    \upClAlg(\dsd)=\upClAlg(\rsd)[x_{-i}^\pm]_{i\in [1,r]}. 
\end{align*}
The restriction of $\dkappa$
to the subalgebra $(\dkappa)^{-1}(R_q[G^{w,u}/H])$ is our previous anti-isomorphism
\begin{align*}
    \rkappa:\upClAlg(\rsd)\simeq R_q[G^{w,u}/H].
\end{align*}

\begin{rem}\label{rem:homogeneous-y-variable}
We have $y_k(\dsd)=y_k(\rsd)$, or, equivalently, $y_k(\dsd)$ is homogeneous of degree in $P\times \{0\}$. To see this, recall that the quantum cluster algebra structure  
    $\rkappa: \bUpClAlg(\rsd)\simeq R_q[G^{w,u}/H] \simeq (S^+_w \bowtie S^-_u)$, constructed in 
    \cite[Theorem 6.2]{goodearl2016berenstein}, uses the symmetric CGL extension presentation of $(S^+_w \bowtie S^-_u)$ with respect to the action of the torus $H$. As a result, all cluster variables are homogeneous with respect to this action, i.e., the $P \times \{0\}$ grading of $R_q[G^{w,u}/H]$.
\end{rem}

Recall that $\uzeta$ and $\ueta$ are reduced words for $w$ and $u$, respectively. Choose the unshuffled signed word $\ubj:=(\uzeta,-\ueta)$ and construct the seed $\dsd(\ubj)$.

Berenstein and Zelevinsky constructed a quantum seed on $R_q[G^{u,w}]$, denoted $\sd^{\BZ}$. Its exchange matrix $\tB^{\BZ}$ is an $\ddI\times \ddI_{\ufv}$-matrix, such that $\tB^{\BZ}=-\tB(\dsd(\ubj))$, see \cite[Eq. (8.7)]{BerensteinZelevinsky05} and Eq. \eqref{eq:signed-word-B-matrix}. Its cluster variables are realized as the generalized quantum minors:
\begin{align}\label{eq:BZ-cluster-var}
    \begin{split}
        x_{-i}(\sd^{\BZ})&=\Delta_{\varpi_i,w^{-1}\varpi_i},\qquad i\in [1,r],\\
    x_{k}(\sd^{\BZ})&=\Delta_{\varpi_{\zeta_k},(w_{>k})^{-1}\varpi_{\zeta_k}},\qquad k\in [1,l(w)],\\
     x_{k+l(w)}(\sd^{\BZ})&=\Delta_{u_{\leq k}\varpi_{\eta_k},\varpi_{\eta _k}},\qquad k\in [1,l(u)].
    \end{split}
\end{align} 
Here, $x_{-i}(\sd^{\BZ})$ is frozen for each $i\in [1,r]$. Moreover, if $k\in [1,r]$ satisfies $k=k^{\max}$ in $[1,l(w)+l(u)]$ with respect to $\ubj$ (see Subsection \ref{sec:ca-sw}), then $x_{k}(\sd^{\BZ})$ is frozen and equals $\Delta_{\varpi_{\zeta_k},(w)^{-1}\varpi_{\zeta_k}}$; if $k\in [1,l(u)]$ satisfies $k+l(w)=(k+l(w))^{\max}$ in $[1,l(w)+l(u)]$ with respect to $\ubj$, then $x_{k+l(w)}(\sd^{\BZ})$ is frozen and equals $\Delta_{u\varpi_{\eta_k},\varpi_{\eta _k}}$. Consequently, the frozen variables are independent of the choice of the reduced words $\uzeta,\ueta$.

Denote the $P\times P$-grading of $x_{k}(\sd^{\BZ})$ by $(\delta^*_k,\delta_k)$, $\forall k\in \ddI$. $\Lambda(\sd^{\BZ})$ is given by the following:
\begin{align*}
    \Lambda^{\BZ}_{kj}=\langle \delta^*_k,\delta^*_j\rangle -\langle \delta_k,\delta_j\rangle,\ \forall k>j.
\end{align*}

Let us make $\dsd(\ubj)$ into a quantum seed such that $\Lambda(\dsd(\ubj))=-\Lambda(\sd^{\BZ})$. Then we have an anti-isomorphism of $\K$-algebras
\begin{align*}
    \dkappa^{\BZ}:\upClAlg(\dsd(\ubj)) \simeq R_q[G^{u,w}]
\end{align*}
such that, $\dkappa x_k(\dsd(\ubj))=x_k(\sd^{\BZ})$, 
$\forall k\in \ddI$. Note that we can identify $\sd^{\BZ}$ with $\dsd(\ubj)^\op$. So we have the following $\K$-algebra isomorphism 
\begin{align*}
    \dkappa^{\BZ}\iota^\op:\upClAlg(\dsd(\ubj)^\op) \simeq R_q[G^{u,w}].
\end{align*}

By \cite{BerensteinZelevinsky05}, the dual Lusztig involution $\overline{(\ )}$ on $R_q[G^{u,w}]$ keeps the generalized quantum minors invariant. So it is identified with the bar involution on $\upClAlg(\dsd^\op)$ via $\dkappa\iota^\op$ and identified with the bar involution on $\upClAlg(\dsd)$ via $\dkappa$. An example can be found in \cite[Example 2.9]{BerensteinZelevinsky05}.

The map $\dkappa$ provides a $P\times P$-grading on $\upClAlg(\dsd(\ubj))$. By \cite[Section 9.2]{goodearl2016berenstein}, we can realize $\upClAlg(\rsd(\ubj))$ as the subalgebra of $\upClAlg(\dsd(\ubj))$ spanned by the elements having degrees in $P \times\{0\}$. In this realization, we have $x_k(\rsd(\ubj))=x_k(\dsd(\ubj))\cdot \factor^{\BZ}_k$, where $k\in \dI$ and $\factor^{\BZ}_k$ is a Laurent monomial in $x_{-i}(\dsd(\ubj))$, $i\in [1,r]$, uniquely determined by $x_k(\dsd(\ubj))$. Then $\dkappa$ restricts to a $\K$-algebra anti-isomorphism
\begin{align}\label{eq:BZ-seed}
    \rkappa^{\BZ}:\upClAlg(\rsd(\ubj))\simeq R_q[G^{u,w}/H].
\end{align}
An example could be found in \cite[Example 9.2]{goodearl2016berenstein}.

An anti-isomorphism 
\begin{equation}
\label{eq:Psi}
\Psi:R_q[G^{u,w}]\rightarrow R_q[G^{w,u}]
\end{equation}
was constructed in \cite{goodearl2016berenstein}. It is given in the form of a composition $\zeta_{w^{-1},u^{-1}}\circ (h_*\cdot) (\omega S)$. Here 
\[
\zeta_{w^{-1},u^{-1}} : R_q[G^{w^{-1},u^{-1}}] \rightarrow R_q[G^{w,u}]
\]
is an anti-isomorphism which is a quantization of the twist map of Fomin--Zelevinsky \cite{fomin1999double}. The quantum twist map was constructed in Section 8 of \cite{goodearl2016berenstein}, see Theorem 8.4 therein for its properties. The isomorphisms 
\[
(\omega S) : R_q[G^{u,w}] \stackrel{\simeq}{\rightarrow} R_q[G^{w^{-1},u^{-1}}],\ (h_*\cdot):R_q[G^{w^{-1},u^{-1}}] \stackrel{\simeq}{\rightarrow}R_q[G^{w^{-1},u^{-1}}],
\]
were constructed in Propositions 7.2 and 7.3 in \cite{goodearl2016berenstein},
where their properties were also described. All of the maps described above are independent of the choice of reduced words for $w$ and $u$. 

In terms of initial seeds, 
the map $\dkappa:\upClAlg(\dsd)\simeq R_q[G^{w,u}]$ is defined by
\begin{align}\label{eq:dkappa-def}
    \dkappa(x_k(\dsd))=\Psi x_k(\sd^{\BZ}).
\end{align}
Then our previous map $\rkappa:\upClAlg(\rsd)\simeq R_q[G^{w,u}/H]$ is determined $\rkappa(x_k(\rsd))=\Psi (x_k(\sd^{\BZ})\cdot \factor^{\BZ}_k)$. See \cite[Theorem 8.4, Proposition 9.4]{goodearl2016berenstein}. 

\subsection{Unique Kazhdan-Lusztig type bases}

Choose $\psi=0$ in the definition of $\sigma'$ and $\cor$ on $U^-_q[w]$, see \eqref{eq:cor}, that is, $\cor(x)=q^{-\frac{1}{4}\langle \wt x,\wt x\rangle}x$ for homogeneous $x$. Then $\varphi^+_w \rkappa W_k=\cor \dPBW_{\beta_{d(k)}}$, $\forall k\in [1,l(w)]$.

Recall that for every $\uc\in \N^{[1,l(w)]}$, we have the standard monomial $\std(\uc)\in\bUpClAlg(\rsd)$. Denote $\bStd(\uc):=\overline{(\std(\uc))}$. Using the anti-isomorphism $\varphi^+_w \rkappa$ and \eqref{eq:cor-dual-PBW}, we obtain the following.
\begin{align}\label{eq:std-dual-PBW}
    \begin{split}
        \varphi^+_w \rkappa \bStd(\uc)&= \varphi^+_w \rkappa \left(q^{-\Hf\sum_{k>j}\Lambda(\rsd)(c_k \deg W_k, c_j \deg W_j)} W_{l(w)}^{c_{l(w)}}\cdots W_1^{c_1}\right)\\
        &= q^{\Hf\sum_{1 \leq j < k \leq l(w)}\nu(c_k\wt W_k,c_j\wt W_j)}(\cor \dPBW_{\beta_{d(1)}})^{c_1}\cdots \cor (\dPBW_{\beta_{d(l(w))}})^{c_{l(w)}}\\
        &=q^{-\Hf\sum_{1 \leq j < k \leq l(w)}\langle c_k\beta_{k}, c_j\beta_{j}\rangle }(\cor \dPBW_{\beta_{d(1)}})^{c_1}\cdots \cor (\dPBW_{\beta_{d(l(w))}})^{c_{l(w)}}\\
        &=\cor (\dPBW(\uc)).       
    \end{split}
\end{align}

Let $\tTri$ denote the common triangular basis of $\bUpClAlg(\rsd)$. By Theorem \ref{thm:bases_dBS}, $\tTri$ is $(<_{\lex},q^{-\Hf}\Z[q^{-\Hf}])$-unitriangular to $\std$. Then it is $(<_{\lex},q^{\Hf}\Z[q^{\Hf}])$-unitriangular to $\bStd:=\{\bStd(\uc)|\uc\in \N^{[1,l]}\}$. Note that, for $\uc^+\in \N^{[1,l(w)]}$ and $\uc^-\in \N^{[l(w)+1,l]}$, we have $\bStd(\uc^++\uc^-)=[\bStd(\uc^-)*\bStd(\uc^+)]$. 

Denote $\tTri^+:=\{\tTri(\uc^+)|\uc^+\in \N^{[1,l(w)]}\}$ and $\tTri^-:=\{\tTri(\uc^-)|\uc^-\in \N^{[l(w)+1,l]}\}$. Then we have $\tTri^\pm(\uc^\pm)\in \bStd(\uc^\pm)+\sum_{ \ud^\pm <_{\lex}\uc^\pm }q^{\Hf}\Z[q^{\Hf}] \bStd(\ud)$, for any $\uc^+\in \N^{[1,l(w)]}$ and $\uc^-\in \N^{[l(w)+1,l]}$. We deduce the following $(<_{\lex},q^{\Hf}\Z[q^{\Hf}])$-unitriangular decomposition:
\begin{align}\label{eq:split-tri-to-std}
    [\tTri^-(\uc^-)*\tTri^+(\uc^+)]&\in \bStd(\uc^+ + \uc^-)+\sum_{\ud<_{\lex}( \uc^+ + \uc^-)}q^{\Hf}\Z[q^{\Hf}]\bStd(\ud).
\end{align}
In particular, $\tTri^{-,+}:=\{[b^-*b^+]|b^\pm \in \tTri^\pm\}$ is a basis of $\bUpClAlg(\rsd)$. Now \eqref{eq:split-tri-to-std} implies the following.
\begin{lem}\label{lem:sign-split-KL}
For any $\uc\in \N^{[1,l]}$, $\tTri(\uc)$ has a $(<_{\lex},q^{\Hf}\Z[q^{\Hf}])$-unitriangular decomposition in $\tTri^{-,+}$.
\end{lem}
\begin{lem}\label{lem:tri-to-dual-canonical}
    We have $\varphi^+_w\rkappa \tTri(\uc)=\cor\dCan(\uc)$, $\forall \uc\in \N^{[1,l(w)]}$.
\end{lem}
\begin{proof} Recall that $\tTri(\uc)$ is determined by its bar-invariance and the $(<_{\lex},q^{\Hf}\Z[q^{\Hf}])$-unitriangularity of its decomposition into $\{ \bStd(\ud) \mid \ud<_{\lex} \uc\}$. Similarly, $\dCan(\uc)$ is determined by its bar-invariance and the $(<_{\lex},q^{\Hf}\Z[q^{\Hf}])$-unitriangularity of its decomposition into $\{ \dPBW(\ud) \mid \ud<_{\lex} \uc\}$. In addition, 
    we have $\varphi^+_w\rkappa \bStd(\ud)=\cor \dPBW(\ud)$, see \eqref{eq:std-dual-PBW}. The desired claim follows.
\end{proof}

\begin{lem}\label{lem:half-triangular-independence}
    The subsets $\varphi^+_w\rkappa \tTri^+ \subset U^-_q[w]$ and $\varphi^-_w\rkappa \tTri^-\subset U^+_q[u]$ are independent of the choice of the reduced words for $w$ and $u$. 
\end{lem}
\begin{proof}
    By Lemma \ref{lem:tri-to-dual-canonical}, we have $\varphi^+_w\rkappa \tTri^+=\cor  \dCan$, where $\dCan$ denotes the dual canonical basis of $U^-_q[w]$. By Lemma \ref{lem:dCan-inv}, it is independent of the choice of reduced word for $w$. 

    Now consider the dual canonical basis for $U^-_q[u]$. Denote $i=|\bi_{k}|$ for $k\in [l(w)+1,l]$ as before. Note that $-\omega$ is an isomorphism from $U^+_q(\frg)$ to $U^-_q(\frg)$. By \eqref{eq:om-prop}, we have 
    \begin{align*}
        (-\omega) (X^+_{\beta'_{d(k)}})&=q_{i}(-1)^{\rho^\vee(\beta'_{d(k)})}q^{-\langle \beta'_{d(k)},\rho\rangle}X^-_{\beta'_{d(k)}}.
    \end{align*}
    Define $\cor$ on $U^-_q[u]$ such that $\cor x=q^{-\frac{1}{4}\langle \wt x,\wt x\rangle}x$ for homogeneous $x$. We obtain
    \begin{align*}
        (-\omega) \varphi^-_u \rkappa W_k&=(-1)^{\rho^\vee(\beta'_{d(k)})}q^{-\langle \beta'_{d(k)},\rho\rangle}q_i^{\Hf}(q_i^{-\Hf}-q_i^{\Hf})X^-_{\beta'_{d(k)}}\\
        &=(-1)^{\rho^\vee(\beta'_{d(k)})}q^{-\langle \beta'_{d(k)},\rho\rangle} \cor \dPBW_{\beta'_{d(k)}}.
    \end{align*}  
For $\uc\in \N^{[l(w)+1,l]}$, denote $\gamma=\sum_{k=l(w)+1}^{l} c_{k}\beta'_{d(k)}$ and $\uc'\in \N^{[1,l(u)]}$ such that $c'_k=c_{k+l(w)}$.
Consider the dual PBW basis and dual canonical basis for $U^-_q[u]$. Using similar arguments as before, we have the following.
\begin{align*}
    (-\omega) \varphi^-_u \rkappa \bStd(\uc)&=(-1)^{\rho^\vee(\gamma)}q^{-\langle \gamma,\rho\rangle} \cor \dPBW(\uc'),\\
    (-\omega) \varphi^-_u \rkappa \tTri^-(\uc)&=(-1)^{\rho^\vee(\gamma)}q^{-\langle \gamma,\rho\rangle} \cor \dCan(\uc').
\end{align*}
The choice independence of $\tTri^-$ then follows from that of the dual canonical basis $\dCan$ for $U^-_q[u]$.
\end{proof}

Combining Lemma \ref{lem:sign-split-KL} and Lemma \ref{lem:half-triangular-independence}, we obtain the following result. 

\begin{prop}\label{prop:unique-CGL-basis}
    The basis $\rkappa \tTri$ of $S^+_w \bowtie S^-_u$ is independent of the choice of reduced words $\uzeta$ and $\ueta$.
    \end{prop}

Let $\rTri$ denote the common triangular basis of $\upClAlg(\rsd)$. 
\begin{thm}\label{thm:unique-rdBC-basis}
    The basis $\rkappa \rTri$ of $R_q[G^{w,u}/H]$ is independent of the choice of reduced words for $w$ and $u$.
    \end{thm}
    \begin{proof}
        The images of the frozen variables $\rkappa x_j(\rsd)$, for $j\in \dI_{\fv}$, are independent of the chosen reduced words. Indeed, they are uniquely determined by the corresponding frozen variables $\dkappa x_j(\dsd)$, which are themselves independent of these choices. The desired statement then follows from Proposition \ref{prop:unique-CGL-basis} by taking the localization such that the frozen variables become invertible.

    \end{proof}

Finally, let $\dTri$ denote the common triangular basis of $\upClAlg(\dsd)$. 
\begin{thm}\label{thm:unique-fdBC-basis}
    The basis $\dkappa \dTri$ of $R_q[G^{w,u}]$ is independent of the choice of reduced words for $w$ and $u$.
    \end{thm}
\begin{proof}
    Recall that we have an algebra inclusion $\upClAlg(\rsd)\subset \upClAlg(\dsd)$, which realizes $\upClAlg(\rsd)$ as the subalgebra spanned by the $P\times \{0\}$-graded components of the $P\times P$-graded algebra $\upClAlg(\dsd)$. For any cluster monomial $z$ of $\upClAlg(\dsd)$, there exists a frozen factor $p$, such that the localized cluster monomial $z\cdot p$ belongs to $\upClAlg(\rsd)$. Then it follows from the definition that for any triangular basis element $b\in \rTri$ and frozen factor $p\in \frGp$, we have $p\cdot b\in \dTri$, i.e., we have $\dTri=\{p\cdot b| b\in \rTri, p\in \frGp\}$, see the base change in \cite[Proposition 4.18]{qin2023analogs}.

    Note that the images of the frozen variables $\dkappa x_{j}(\dsd)$, $j\in I_\fv(\dsd)$, are independent of the choice of the reduced words, and $\dkappa\rTri$ is independent of the choice by  Theorem \ref{thm:unique-rdBC-basis}. The desired claim follows.
\end{proof}

\subsection{Canonical cluster structure on quantum double Bruhat cells}

Let $\uzeta'$ and $\ueta'$ be another choice of the reduced words for $w$ and $u$, respectively. Construct $\ubi'=((\uzeta')^\op,\ueta')$. Using the quantization matrix provided by \cite{goodearl2016berenstein}, see \eqref{eq:Lambda-GY}, we obtain quantum seeds $\dsd':=\dsd(\ubi')$, $\rsd':=\rsd(\ubi')$, and anti-isomorphisms $\dkappa':\upClAlg(\dsd')\simeq R_q[G^{w,u}]$, $\rkappa':\upClAlg(\rsd')\simeq R_q[G^{w,u}/H]$.

Recall that there is a permutation-mutation sequence $\seq^\sigma=\seq^\sigma_{\ubi',\ubi}$ associated with braid moves. Denote $\dsd''=\seq^\sigma \dsd$, $\rsd''=\seq^\sigma \rsd$. Recall that we have $\dsd''=\dsd'$ at the classical level, see \cite[Proposition 3.7]{shen2021cluster}. A priori, the quantization matrices $\Lambda(\dsd'')$ and $\Lambda(\dsd')$ might be different. So we have the following diagrams.

\begin{align}\label{eq:dBC-mutation-diagram}
    \begin{array}{ccc}
    \upClAlg(\dsd'') & \stackrel[\text{change}]{\text{quantization}}{\dashrightarrow} & \upClAlg(\dsd')\\
    \simeqd(\seq^{\sigma})^{*} &  & \simeqd\dkappa'\\
    \upClAlg(\dsd)  & \overset{\dkappa}{\xrightarrow{\sim}} & R_q[G^{w,u}]
    \end{array}
    \end{align}

\begin{align}\label{eq:rdBC-mutation-diagram}
    \begin{array}{ccc}
    \upClAlg(\rsd'') & \stackrel[\text{change}]{\text{quantization}}{\dashrightarrow} & \upClAlg(\rsd')\\
    \simeqd(\seq^{\sigma})^{*} &  & \simeqd\rkappa'\\
    \upClAlg(\rsd) & \overset{\rkappa}{\xrightarrow{\sim}} & R_q[G^{w,u}/H]
    \end{array}
    \end{align}

By \cite[Theorem 1.1]{shen2021cluster}, at the classical level, we have $\upClAlg(\dsd'')=\upClAlg(\dsd')$ and $\dkappa'=\dkappa(\seq^\sigma)^*$. Restricting to the subalgebras spanned by the $P\times \{0\}$-graded components, we obtain $\upClAlg(\rsd'')=\upClAlg(\rsd')$ and $\rkappa'=\rkappa(\seq^\sigma)^*$ at the classical level. We aim to lift these results to the quantum level.

\begin{thm}\label{thm:uni-cluster-rdBC}
We have $\rsd'=\rsd''$, i.e., $\upClAlg(\rsd'')=\upClAlg(\rsd')$. Moreover, $\rkappa'=\rkappa(\seq^\sigma)^*$.
\end{thm}
\begin{proof} The proof of this theorem is similar to that of \cite[Theorem 6.4]{qin2024infinite}.

(i) 
Let $W'_{[j,k]}$ and $W''_{[j,k]}$ denote the interval variables of $\upClAlg(\rsd')$ and $\upClAlg(\rsd'')$, respectively. We claim that $\rkappa(\seq^\sigma)^*W''_{[j,k]}=\rkappa' W'_{[j,k]}$.

Note that $(\seq^\sigma)^*W''_{[j,k]}$ are quantum cluster variables of $\upClAlg(\rsd)$. Therefore, each of them belongs to the common triangular basis $\rTri$ of $\upClAlg(\rsd)$. It follows that $\rkappa (\seq^\sigma)^*W''_{[j,k]}$ belongs to the basis $\rkappa \rTri$ of $R_q[G^{w,u}/H]$. On the other hand, since the quantum cluster variable $W'_{[j,k]}$ belongs to the common triangular basis $\rTri'$ of $\bUpClAlg'$, $\rkappa' W'_{[j,k]}$ belongs to the basis $\rkappa' \rTri'$ of $R_q[G^{w,u}/H]$. By Theorem \ref{thm:unique-rdBC-basis}, we have $\rkappa\rTri=\rkappa'\rTri'$. Therefore, $\rkappa^{-1}\rkappa'W'_{[j,k]}$ belongs to $\rTri$. Since the elements of $\rTri$ have distinct degrees in $\LP(\rsd)$, if we could show that $\rkappa^{-1}\rkappa'W'_{[j,k]}$ and $(\seq^\sigma)^*W''_{[j,k]}$ have the same degree, then they must be equal.

It suffices to compare their degrees at the classical level. Take $\kk=\C$ now. We have $\rsd'=\rsd''$, and $\rkappa'=\rkappa (\seq^\sigma)^*$ by \cite[Theorem 1.1]{shen2021cluster}. Therefore, $\rkappa^{-1}\rkappa'W'_{[j,k]}=(\seq^\sigma)^*W''_{[j,k]}$ at the classical level. In particular, they have the same degree. The desired claim is verified.

(ii) Denote $\rkappa'':=\rkappa (\seq^\sigma)^*$. We have shown the $\K$-algebra isomorphism $(\rkappa'')^{-1}\rkappa':\bUpClAlg'\simeq \bUpClAlg''$ sends $W'_{[j,k]}$ to $W''_{[j,k]}$. In particular, the initial cluster variables $x'_k=W'_{[k^{\min},k]}$ are sent to $x''_k=W''_{[k^{\min},k]}$. By comparing $x'_j*x'_k =q^{\Lambda(\rsd')_{jk}}x'_k*x'_j$ and $x''_j*x''_k =q^{\Lambda(\rsd'')_{jk}}x''_k*x''_j$, we deduce that $\Lambda(\rsd')=\Lambda(\rsd'')$. Therefore, $\rsd'=\rsd''$ and thus $\upClAlg(\rsd')=\upClAlg(\rsd'')$. Combining with (i), we deduce that $\rkappa(\seq^\sigma)^*=\rkappa'$.
\end{proof}

\begin{thm}\label{thm:uni-cluster-dBC}
    We have $\dsd'=\dsd''$ as quantum seeds, i.e., $ \upClAlg(\dsd'')= \upClAlg(\dsd')$. Moreover, $\dkappa'=\dkappa(\seq^\sigma)^*$.
    \end{thm}
\begin{proof}
Recall that we have identified $\upClAlg(\rsd)$ as the subalgebra of $\upClAlg(\dsd)$ spanned by the $P\times \{0\}$-graded components. The cluster variable $x_k(\dsd'')$ takes the form $\factor''_k\cdot x_k(\rsd'')$ for $\factor''_k$ a Laurent monomial in $x_{-i}$, $i\in [1,r]$. Similarly, $x_k(\dsd')$ takes the form  $\factor'_k\cdot x_k(\rsd')$ for $\factor'_k$ a Laurent monomial in $x_{-i}$.

(i) Note that $x_{-i}$ are frozen variables. We have $\dkappa(\seq^\sigma)^*x_{-i}=\dkappa'x_{-i}$, see \eqref{eq:dkappa-def}. Then $(\dkappa')^{-1}\dkappa (\seq^\sigma)^* \factor''_k$ is a frozen factor in $\upClAlg(\dsd')$. We claim that it equals $\factor'_k$. 

It suffices to verify the claim at the classical level. Note that, $\forall k\in [1,l]$, we have $\dkappa (\seq^\sigma)^* x_k(\dsd'')=\dkappa' x_k(\dsd')$ at the classical level by \cite[Theorem 1.1]{shen2021cluster}, and we have seen that $\rkappa (\seq^\sigma)^* x_k(\rsd'')=\rkappa' x_k(\rsd')$. We deduce that $\dkappa (\seq^\sigma)^* \factor''_k=\dkappa' \factor'_k$ at the classical level.

(ii) In (i), we have shown $\dkappa (\seq^\sigma)^* \factor''_k=\dkappa' \factor'_k$. In addition, $\rkappa (\seq^\sigma)^* x_k(\rsd'')=x_k(\rsd')$ by Theorem \ref{thm:uni-cluster-rdBC}. We deduce that $\dkappa (\seq^\sigma)^* x_k(\dsd'')= \dkappa' x_k(\dsd')$. Then we deduce that $\dsd''$ and $\dsd'$ have the same quantization matrix, by using the same arguments as in the proof (i) for Theorem \ref{thm:uni-cluster-rdBC}. 

Therefore, $\dsd''=\dsd'$. Combined with $\dkappa (\seq^\sigma)^* x_k(\dsd'')= \dkappa' x_k(\dsd')$ and $\dkappa(\seq^\sigma)^*x_{-i}=\dkappa'x_{-i}$, we obtain $\dkappa (\seq^\sigma)^*= \dkappa'$.
\end{proof}

Since the quantum seed of Berenstein--Zelevinsky on $R_q[G^{u,w}]$ and that of Goodearl--Yakimov on $R_q[G^{w,u}]$ are related by the anti-isomorphism $\Psi$
in \eqref{eq:Psi},
Theorem \ref{thm:uni-cluster-dBC} implies the following:
\begin{thm}\label{thm:uni-cluster-BZ-dBC}
The quantum seeds $\sd^{\BZ}$ and $(\sd^{\BZ})'$ of $R_q[G^{u,w}]$ associated with the signed words $(\uzeta,-\ueta)$ and $(\uzeta',-\ueta')$, respectively, are related by permutations and mutations.
\end{thm}

\subsection{A proof of the full Berenstein--Zelevinsky conjecture}
As before, let $\ueta,\uzeta$ denote two reduced words of Weyl group elements $u,w$, respectively, and $l:=l(u)+l(w)$. Let $\ubi$ denote an arbitrary shuffle of $\uzeta$, $-\ueta$. Recall from the introduction
and Section \ref{sec:ca-sw} that the associated (candidate) quantum seed $\sd^\BZ$ of \cite{BerensteinZelevinsky05} for $R_q[G^{u,w}]$ has the following quantum cluster variables:
\begin{align}\label{eq:BZ-cluster-var-general}
    \begin{split}
        x_{-i}(\sd^{\BZ})&=\Delta_{\varpi_i,w^{-1}\varpi_i},\qquad i\in [1,r],\\
    x_{k}(\sd^{\BZ})&=\Delta_{u_{\leq k}\varpi_{|\bi_k|},w^{-1}_{> k}\varpi_{|\bi_k|}},\qquad k\in [1,l],
    \end{split}
\end{align} 
where 
\[
u_{\leq k}:=s_{-\bi_1}s_{-\bi_2}\cdots s_{-\bi_k} \quad \mbox{and} \quad  
w^{-1}_{> k}:=s_{\bi_l}s_{\bi_{l-1}}\cdots s_{\bi_{k+1}},
\]
and $s_{-a}$ denotes the identity for $a\in [1,r]$.  

Let $\ubj$ denote the unshuffled signed word $(\uzeta,-\ueta)$. We can always choose a finite sequence of flips from $\ubj$ to $\ubi$, such that each flip changes a pair of adjacent letters $(m,-n)$ sitting in positions $[1,l]$ to $(-n,m)$, where $m,n\in [1,r]$. Denote the resulting sequence of signed words by $\ubi_0:=\ubj$, $\ubi_1$, ... ,$\ubi_N=\ubi$ where $N$ is a positive integer. Let $\sd_t$, $t\in [0,N]$ be 
the associated Berenstein--Zelevinsky quantum seeds. 

\begin{prop}
\label{lem:flip-BZ-seeds}
The Berenstein--Zelevinsky quantum seeds $\sd_t$ are connected by permutations and mutations. 
\end{prop}
\begin{proof}
Fix $t\in [0,N-1]$. For some $k\in [1,l-1]$ and $m,n\in [1,r]$, we have $(\ubi_t)_{[k,k+1]}=(m,-n)$ and it is flipped to $(\ubi_{t+1})_{[k,k+1]}=(-n,m)$. 

When $m \neq n$, using \eqref{eq:BZ-cluster-var-general}, we see that the quantum seeds $\sd_t$ and $\sd_{t+1}$ share the same collection of cluster variables, such that $x_{i}(\sd_{t+1})=x_{\sigma(i)}(\sd_{t})$ for the transposition 
$\sigma= (k \, k+1)$. One easily verifies that the corresponding exchange matrices (equal to minus those in \eqref{eq:signed-word-B-matrix}) also match under the transposition $\sigma$, and so the two seeds are related by a permutation.

When $m=n$, using \eqref{eq:BZ-cluster-var-general}, we see that the quantum seeds $\sd_t$ and $\sd_{t+1}$ satisfy $x_{i}(\sd_{t+1})=x_{i}(\sd_{t})$ for $i\in \ddI\backslash \{k\}$, and 
\[
x_{k}(\sd_{t})=
\Delta_{u_{\leq k}\varpi_{m},w^{-1}_{> k}\varpi_{m}}, 
\quad 
x_{k}(\sd_{t+1})=\Delta_{u_{\leq k}s_m\varpi_{m},w^{-1}_{> k}s_m\varpi_{m}}.
\]
Moreover, since $\ueta$, $\uzeta$ are reduced words, we have $l(u_{\leq k}s_m)=l(u_{\leq k})+1$, $l(w^{-1}_{> k}s_m)=l(w^{-1}_{> k})+1$. By \cite[Proposition 3.2]{GeissLeclercSchroeer11}, we have the following relation
\begin{multline}\label{eq:GLS-mutation}
    \Delta_{u_{\leq k}s_m\varpi_m,w^{-1}_{> k}s_m\varpi_m}\Delta_{u_{\leq k}\varpi_m,w^{-1}_{> k}\varpi_m}
    \\
    =q^{-1} \Delta_{u_{\leq k}s_m\varpi_m,w^{-1}_{> k}\varpi_m}\Delta_{u_{\leq k}\varpi_m,w^{-1}_{> k}s_m\varpi_m}
    +\prod_{p \in [1,r] \backslash \{m \}}\Delta_{u_{\leq k}\varpi_p,w^{-1}_{> k}\varpi_p}^{-c_{pm}}.
\end{multline}
(The paper \cite{GeissLeclercSchroeer11} considers symmetric Kac--Moody algebras $\g$ but the proof of this particular result does not use the symmetric assumption.)
In this case we have $k[1]=k+1$. Recall that the corresponding Berenstein--Zelevinsky exchange matrix is $B^{\BZ}(\sd_t)=-B(\dsd(\ubi_t))$, where the exchange matrix $B(\dsd(\ubi_t))$ is given by \eqref{eq:signed-word-B-matrix}. Its entries in the $k$-th column are
\begin{align}
    \begin{split}        
b_{jk}(\sd_t)&=-1,\quad \text{if $j=k[1]=  k+1$ or $j= k[-1]$},\\
b_{jk}(\sd_t)&=-c_{pm },\quad \text{if $|\bi_j|=p\neq m$ and 
$j<k<k[1]=k+1<j[1]$},\\
b_{jk}(\sd_t)&=0,\quad \text{for all other $j$}.
    \end{split}
\end{align}
Note that for every $p \in [1,r]\backslash\{m\}$ there exists 
a unique $j \in [-r,-1] \sqcup [1,k-1]$ with the property in the second equality in the above display because 
$\bi_j = -p$ for some $j \in [-r,-1]$. Note that 
in that equality $j[1]$ may equal $\infty$. 
Combined with \eqref{eq:BZ-cluster-var-general}, we can translate \eqref{eq:GLS-mutation} into the following equality.
\begin{align*}
    x_k(\sd_{t+1})x_k(\sd_{t})=q^{-1}x_{k[1]}(\sd_{t})x_{k[-1]}(\sd_{t})+\prod_{j\neq k}x_{j}(\sd_{t})^{[-b_{jk}]_+}.
\end{align*}
On the other hand, applying the mutation $\mu_k$ to $\sd_t$, we obtain
\begin{align*}
    x_k(\mu_k\sd_{t})x_k(\sd_{t})=q^{\alpha}x_{k[1]}(\sd_{t})x_{k[-1]}(\sd_{t})+q^{\beta}\prod_{j\neq k}x_{j}(\sd_{t})^{[-b_{jk}]_+}.
\end{align*}
for some $\alpha,\beta\in \frac{\Z}{2}$. Since both $x_k(\mu_k\sd_{t})$ and $x_k(\sd_{t+1})$ are bar-invariant in $\LP(\sd_t)$, we must have $\alpha=-1$ and $\beta=0$. It follows that $x_k(\mu_k\sd_{t})=x_k(\sd_{t+1})$ and, consequently, $\mu_k\sd_{t}=\sd_{t+1}$.
\end{proof}

\noindent
{\em{Proof of Theorem A.}}
Combining Theorem \ref{thm:uni-cluster-BZ-dBC} and Proposition \ref{lem:flip-BZ-seeds}, we obtain that all Berenstein-Zelevinsky seeds for the pair $(u,w)$ associated to arbitrary signed words are connected by permutations and mutations. Since by \cite{goodearl2016berenstein}, the 
seeds for unshuffled signed words are quantum seeds of $R_q[G^{u,w}]$, the same is true for all seeds. By 
\cite{goodearl2016berenstein} we have 
\[
R_q[G^{u,w}] = \upClAlg(\sd^{\BZ}(\ubi)) = \clAlg(\sd^{\BZ}(\ubi))
\]
for unshuffled words, and applying again 
Theorem \ref{thm:uni-cluster-BZ-dBC} and Proposition \ref{lem:flip-BZ-seeds} gives that the same holds for all signed words. 
\qed
\section{An approach to quantum partial compactifications}	\label{sec:key}
In this section we prove Theorem C.
\subsection{Intersection of localizations of domains}

Let $A$ be a left or right Ore domain (i.e., $0$ is the only zero divisor and $A$ has a left or right quotient ring $Q(A)$, see e.g. 
\cite[\S 2.1.14]{McConnell-Robson}. The ring $A$ can be noncommutative). 
For example, every left Noetherian domain is a left Ore domain \cite[Theorem 2.1.15]{McConnell-Robson}.

\begin{definition} 
\label{def:prime}
A {\em{prime element}} of $A$ is a nonzero, nonunit, normal element $\Delta \in A$ (recall Definition \ref{def:normal}) such that $A/(\Delta)$ is a domain, where $(\Delta):= A \Delta = \Delta A$ denotes the principal ideal generated by $\Delta$, i.e., $(\Delta)$ is a completely prime ideal of $A$.  
\end{definition}

For any $a\in A$ and a normal element $\Delta \in A$, we say $\Delta | a$ if $a\in (\Delta)$. The condition that for a prime element $\Delta \in A$, $A/(\Delta)$ is a domain is equivalent to the classical property
\[
\Delta | ab \quad \Rightarrow \quad \Delta | a \,\,\,\,
\mbox{or} \,\,\,\, \Delta | b \quad \mbox{for} \,\, a, b \in A.
\]

\begin{thm}\label{thm:intersection}
Let $J$ denote an index set. Assume that $\Delta_j$, $j\in J$, are prime elements and $E\subset A$ is a left or right denominator set 
\cite[\S 2.1.13]{McConnell-Robson}
such that $\Delta_j \nmid x$ for any $j\in J$ and $x\in E$. Then we have
\begin{align}
\label{eq:inter}
    A[\Delta_j^{-1}]_{j\in J}\cap A[E^{-1}]=A.
\end{align}
\end{thm}
We note that the multiplicative subset of $A$ generated by $\{\Delta_j \mid j\in J\}$ is a left and right denominator set because the elements $\Delta_j \in A$ are normal. Since $A$ is assumed to be a domain, $A$ canonically embeds in $A[\Delta_j^{-1}]_{j\in J}$ and $ A[E^{-1}]$. The assumption that $A$ is a left or right Ore domain gives that $A[\Delta_j^{-1}]_{j\in J}$ and $ A[E^{-1}]$
are canonically subrings of $Q(A)$, and this is 
where the intersection \eqref{eq:inter} takes place. 

In other words, Theorem \ref{thm:intersection} states that for any algebra $B$ such that 
 \[
 A\subset B\subset A[\Delta_j^{-1}]_{j\in J}\cap A[E^{-1}],
 \]
 we have $A=B$.

\begin{proof}
We only treat the case where $E$ is a left denominator set. The proof for a right denominator set is similar.

Any $0\neq b\in A[\Delta_j^{-1}]_{j\in J}\cap A[E^{-1}]$ can be written as
\begin{align*}
    b=a\Delta^{-1}=x^{-1}a'
\end{align*}
where $\Delta$ is a product of $r$ $\Delta_j$'s with possible repetitions, $x\in E$, and $a,a'\in A$. Set $\Delta = \Delta' \Delta_i$, 
where $\Delta'$ is a product of $(r-1)$ $\Delta_j$'s. It follows that
\begin{align*}
    x a= a' \Delta' \Delta_i.
\end{align*}
We deduce that $\Delta_i| a$ from the facts that $\Delta_i$ is prime and $\Delta_i\nmid x$. We can then replace $a$ by $a \Delta_i^{-1} \in A$ and $\Delta$ by $\Delta'$ since $A$ is a domain. Arguing by induction on $r$, we obtain that $a=a_0 \Delta$ for some $a_0\in A$. Consequently, we have
\begin{align*}
    b=a\Delta^{-1}=a_0 \Delta \Delta^{-1}=a_0\in A.
\end{align*}
\end{proof}

\subsection{Applications to cluster algebras}

Choose any initial seed $\sd$. Recall that  $\bLP$ denotes the associated quantum mixed torus-affine algebra whose frozen variables $x_j$, $j\in I_\fv$, are not inverted.
Let $E$ denote the set consisting of the monomials of the initial unfrozen variables. 

Since the quantum torus $\LP$ is a domain, its subalgebras $\bClAlg$ and $\bUpClAlg$ are domains as well. It is obvious from the definitions that any frozen variable $x_j$ is a normal element of the partially compactified quantum cluster algebra $\bClAlg$ and upper quantum cluster algebra $\bUpClAlg$. 

The following is a quantum version of \cite[Proposition 4.3(ii)]{MNTY}.

\begin{lem}
\label{lem:primeU}
Every frozen variable $x_j$ is a prime element of $\bUpClAlg$.
\end{lem}
\begin{proof} Assume that there exist $a,b,c\in \bUpClAlg$, such that $x_j\nmid a,b$ but $ab=x_j c$. View $a,b,c$ as elements in $\bLP$ and write them as rational functions in reduced forms, whose numerators and denominators do not have common divisors. By \eqref{eq:order-bUpClAlg}, their denominators are monomials of the unfrozen variables. 

Note that $x_j^{-1}a,x_j^{-1}b\in \upClAlg$. In addition, we have $\bUpClAlg=\{z\in\upClAlg|\nu_j(z)\geq 0\}$, again by \eqref{eq:order-bUpClAlg}. Then $x_j\nmid a,b$ in $\bUpClAlg$ implies $\nu_j(a)=\nu_j(b)=0$, i.e., the numerators of $a,b$ are not divisible by $x_j$ as polynomials. 

We can choose a monomial $e$ in the unfrozen variables such that $abe$ and $ce$ are polynomials in the initial cluster variables. We get $abe=x_j ce$. But this is impossible since the numerators of $a,b$ are not divisible by $x_j$.
\end{proof}

\begin{theorem}\label{thm:denominator-set}
    Let $A$ be a subalgebra of $\bLP$ which contains all initial cluster variables and $E$ be the multiplicative subset of $A$ generated by all initial unfrozen cluster variables.  
    Then $E$ is a left and right denominator set of $A$, 
    \begin{equation}
    \label{eq:AbLP}
    A[E^{-1}] = \bLP   
    \end{equation}
    and
    \begin{equation}
    x_j \nmid u \quad for \quad j \in I_f, u \in E.
    \label{eq:not-div}
    \end{equation}
    If $x_j \in A$ are prime elements for all $j \in I_\fv$, then 
    \begin{equation}
        \label{eq:Aint}
        A = A[x_j^{-1}]_{j \in I_\fv} \cap \bLP.
    \end{equation}
\end{theorem}
\begin{rem}
    Since the cluster variables in a quantum seed $q$-commute, $E$ consists of scalar multiples of monomials in the unfrozen initial cluster variables. 
\end{rem}
\begin{proof}
We have that $A \subset \bLP$ and every element of $E$ is invertible in $\bLP$. Furthermore, every element of $\bLP$ is a linear combination of left (and right) fractions with numerators equal to monomials in the initial cluster variables (elements of $A$) and denominators equal to monomials in the initial unfrozen cluster variables (elements of $E$). Therefore, $E$ is a left and right denominator set of $A$ and $A[E^{-1}] = \bLP$. 

Since $x_j^{-1} u \notin \bLP$, we have \eqref{eq:not-div}. The last equality \eqref{eq:Aint} now follows from Theorem \ref{thm:intersection}.
\end{proof}

\begin{corollary} 
\label{cor:Uint}
For every upper quantum cluster algebra 
$\upClAlg$,
     we have 
\[
     \bUpClAlg=\bUpClAlg[x_j^{-1}]_{j\in I_\fv}\cap \bUpClAlg[E^{-1}] = \upClAlg \cap \bLP,
\]
where $E$ is the multiplicative subset of $\upClAlg$ generated by all unfrozen initial cluster variables.
\end{corollary}
We offer two proofs, one based on Theorem \ref{thm:denominator-set} and the other on Formula \eqref{eq:order-bUpClAlg}. 
\medskip

\noindent
{\em{First Proof.}} It is easy to see that
\begin{equation}
\label{eq:eq1}
\bUpClAlg[x_j^{-1}]_{j\in I_\fv} = \upClAlg.
\end{equation}
Both equalities now follow from Lemma \ref{lem:primeU} and Theorem \ref{thm:denominator-set}. 
\qed
\medskip

\noindent
{\em{Second Proof.}} Applying Formula \eqref{eq:order-bUpClAlg}, we obtain 
\[
\bUpClAlg = \upClAlg \cap \bLP.
\]
The first equality in the corollary follows from Eq. \eqref{eq:AbLP} applied to $A = \bUpClAlg$ and Eq. \eqref{eq:eq1} in the first proof of the corollary.
\qed

\noindent
{\em{Proof of Theorem C.}} Theorem \ref{thm:denominator-set} implies that 
\[
A = A[x_j^{-1}]_{j \in I_\fv} \cap \bLP.
\]
By the assumption \eqref{eq:loc} in Theorem C, we have $A[x_j^{-1}, j \in I_\fv] = \upClAlg$. So, 
\[
A = \upClAlg \cap \bLP. 
\]
However, by Corollary \ref{cor:Uint},
\[
     \bUpClAlg= \upClAlg \cap \bLP,
\]
so $A = \bUpClAlg$.
\qed

We expect the following conjecture to hold for a broad class of well-behaved cluster algebras, including those whose frozen vertices can be optimized. It does not, however, hold in complete generality; see Remark \ref{rem:failure-conj} for a counterexample. 
\begin{conj}\label{conj:prime_bA}
Every frozen cluster variable $x_j$ is a prime element of $\bClAlg$.
\end{conj}
It has the following important application.
\begin{prop} 
\label{prop:AeqU}
If a cluster algebra $\clAlg$ satisfies Conjecture \ref{conj:prime_bA}, then
\[
\clAlg=\upClAlg \quad \Leftrightarrow \quad \bClAlg=\bUpClAlg. 
\]   
\end{prop}
\begin{proof} 
The left implication is obvious since 
\[
\bClAlg[x_j^{-1}]_{j\in I_\fv} = \clAlg 
\quad \mbox{and} \quad
\bUpClAlg[x_j^{-1}]_{j\in I_\fv} = \upClAlg. 
\]
Next we show the right implication.
Theorem \ref{thm:denominator-set} and the condition in Conjecture \ref{conj:prime_bA} imply  
\[
\bClAlg=\bClAlg[x_j^{-1}]_{j\in I_\fv}\cap \bLP.
\]
Taking again into account that $\bClAlg[x_j^{-1}]_{j\in I_\fv} = \clAlg$, we get that 
\[
\bClAlg = \clAlg \cap \bLP.
\]
The right implication now follows from Eq. \eqref{eq:order-bUpClAlg} 
(or Corollary \ref{cor:Uint}). 
\end{proof} 
\begin{rem}\label{rem:failure-conj} Bucher, Machacek, and Shapiro \cite{BMS} proved that for an acyclic seed with optimized frozen vertices, $\bClAlg = \bUpClAlg$ at the classical level. They also constructed an example with  
$\clAlg = \upClAlg$ and 
$\bClAlg \neq  \bUpClAlg$. It follows from 
Proposition \ref{prop:AeqU} that their partially compactified 
cluster algebra $\bClAlg$ has a frozen variable that is not a prime element. 
\end{rem}
\section{Quantum cluster structure on quantized coordinate ring}\label{sec:cluster-G}
In this section we prove Theorem B.
We work with quantum cluster algebras and the quantum function algebras $R_q[G]$ defined over $\K=\C(q^\Hf)$.

Let $\uzeta,\ueta$ denote two reduced words for the longest Weyl group element $w_0$. Choose the signed word $\ubi=(\uzeta,-\ueta)$ and take the associated quantum seed $\sd=\sd^{\BZ}(\ubi)$ of Berenstein--Zelevinsky. Recall that $\sd^{\BZ}(\ubi)=\dsd(\ubi)^\op$. For simplicity of the notation, when the choice of the signed word is not relevant, we will suppress the argument $\ubi$ when denoting the quantum seed.
We have the isomorphism
\begin{align*}
    \dkappa^{\BZ}\iota^\op:\upClAlg(\sd)\stackrel{\simeq}\longrightarrow R_q[G^{w_0,w_0}].
\end{align*}
From now on, let us identify $\upClAlg(\sd)$ with $R_q[G^{w_0,w_0}]$ via $\dkappa^{\BZ}\iota^\op$. Then the initial cluster variables are generalized quantum minors. The frozen variables are 
\begin{equation}
    \label{eq:frozen}
\Delta_{w_0 \varpi_i, \varpi_i}, \, \, 
\Delta_{\varpi_i, w_0 \varpi_i}, \quad i \in [1,r]. 
\end{equation}
They are the cluster variables $x_j$ indexed by $\{-i \mid i\in [1,r]\}\sqcup \{i^{\max}\mid i\in [1,r]\}$.
\begin{lem} 
\label{lem:frozen-minor-prime}
The frozen variables in \eqref{eq:frozen} are prime elements of $R_q[G]$.
\end{lem}
\begin{proof} The elements in \eqref{eq:frozen} are normal elements of $R_q[G]$ by \cite[Theorem I.8.15]{joseph2012quantum}, \cite[Lemma 4.2]{HLT}. In \cite[Th\'eor\`eme 3]{Joseph-cr} it was proved that 
the principal ideals of $R_q[G]$ generated by them are the ideals denoted by $Q_{s_i w_0, w_0}$ and 
$Q_{w_0, s_i w_0}$ in \cite{joseph2012quantum}. Those ideals were shown to be completely prime in 
\cite[Proposition 8.9]{joseph2012quantum}; this also appears in  \cite[Theorem 4.4]{HLT}, phrased in a slightly different way.
\end{proof} 

By Theorem \ref{thm:uni-cluster-BZ-dBC}, different choices of the reduced words $\uzeta$, $\ueta$, provide different quantum seeds for the same cluster algebra $R_q[G^{w_0,w_0}]=\upClAlg(\sd)$. Then, for any frozen vertex $j$, we can construct a new quantum seed $\sd_j$ of $R_q[G^{w_0,w_0}]$ as follows.
\begin{itemize}
\item    If $j=-i$, we choose a reduced word $\uzeta^{(i)}$ of $w_0$ whose first letter is $i$. Denote $\ubi^{(j)}:=(\uzeta^{(i)},-\ueta)$ and $\sd_j:=\dsd(\ubi^{(j)})^\op$. Note that $b_{jk}(\sd_j)=0$ whenever $k\in I_\ufv\backslash \{j[1]\}$.

Define $\ubi':=\ubi^{(j)}_{[2,l]}$ and $\sd'_j=\dsd(\ubi')^\op$. Define the codimension-$1$ double Bruhat cell $V_j:=G^{{w_0,s_i w_0}}$.

\item If $j=i^{\max}$, we choose a reduced word $\ueta^{(i)}$ of $w_0$ whose last letter is $i$. Denote $\ubi^{(j)}:=(\uzeta,-\ueta^{(i)})$ and $\sd_j:=\dsd(\ubi^{(j)})^\op$. Note that $b_{jk}(\sd_j)=0$ whenever $k\in I_\ufv\backslash \{j[-1]\}$.

Define $\ubi':=\ubi^{(j)}_{[1,l-1]}$ and $\sd'_j=\dsd(\ubi')^\op$. Define the codimension-$1$ double Bruhat cell $V_j:=G^{{w_0 s_i,w_0}}$.
\end{itemize}
Note that all cluster variables $x_k(\sd_j)$ are generalized quantum minors, see \eqref{eq:BZ-cluster-var} and Theorem \ref{thm:uni-cluster-BZ-dBC}. 

Recall that $\nu_j$ denotes the order of vanishing at $x_j=0$ of an element of the upper cluster algebra structure on $R_q[G^{w_0,w_0}]$ for $j \in I_\fv$.
\begin{lem}\label{lem:nu-j-positive}
For any $f\in R_q[G]$, we have $\nu_j(f)\geq 0$ for $j \in I_\fv$.
\end{lem}
\begin{proof}
Let us write the reduced expression of $f$ in $\LP(\sd_j)$ as $f=x^{-d}*\sum_M c_M M$, where $x^d$ and $M$ are cluster monomials of $\sd_j$, and $c_M\in \K$, such that $d\in \N^{\ddI}$ is minimal. So we have $x^d *f=\sum_M c_M M$. 

Assume $\nu_j(f)<0$, then we have $d_j=-\nu_j(f)> 0$, and $\sum c_M M$ is not divisible by $x_j$ as a polynomial in the cluster variables. Then $\sum_M c'_M M$ is nonzero, where $c'_M=c_M$ if $M$ is not divisible by $x_j$, and $c'_M=0$ otherwise.

We have $x^d *f=\sum_M c_M M$ in $R_q[G]$. Therefore, the image of $\sum_M c'_M M$ in $R_q[G]/(x_j)$ must be equal to $0$. Note that $x_j$ is prime in $R_q[G]$, and the quotient $R_q[G]/(x_j)$ is a domain, see Lemma \ref{lem:frozen-minor-prime}.

Moreover, recall that $R_q[V_j]$ is a localization of the quotient $R_q[G]/(x_j)$, see Section \ref{sec:d-Bruhat} and \eqref{eq:BZ-cluster-var}. Then the identification $\upClAlg(\sd_j)=R_q[G^{w_0,w_0}]$ induces the identification $\upClAlg(\sd'_j)=R_q[V_j]$, such that the cluster variables $x_k(\sd'_j)$ are the images of the cluster variables $x_k(\sd_j)$ in the quotient $R_q[G]/(x_j)$. It follows that the image of $\sum_M c'_M M$ is non-zero in $\upClAlg(\sd'_j)$. This contradicts the fact that 
$x^d * f=\sum_M c_M M$ should vanish in $R_q[G]/(x_j)$ because of the assumption that $d_j = - \nu_j(f) >0$. 
Hence, $\nu_j(f)$ must be nonnegative for all 
$j \in I_\fv$.
\end{proof}

\begin{prop}\label{prop:Guv_in_bU}
    We have $R_q[G]\subset \bUpClAlg$.
\end{prop}
\begin{proof}
    Recall that $R_q[G]\subset R_q[G^{w_0,w_0}]=\upClAlg$. Moreover, $\bUpClAlg=\{z\in \upClAlg|\nu_j(z)\geq 0,\ \forall j\in I_\fv\}$, see \eqref{eq:order-bUpClAlg}. The claim follows from Lemma \ref{lem:nu-j-positive}.
\end{proof}

Recall that $R_q[G^{w_0,w_0}]=\clAlg(\sd)=\upClAlg(\sd)$.
\begin{thm}\label{thm:quantized_coordinate_ring}
    We have $R_q[G]=\bUpClAlg(\sd)$.
\end{thm}
\begin{proof}
    Recall that $R_q[G]$ contains all initial cluster variables, The initial frozen cluster variables are given by \eqref{eq:frozen} and they are prime elements of $R_q[G]$ by Lemma \ref{lem:frozen-minor-prime}. Moreover, 
    \[
    R_q[G][\Delta_{w_0 \varpi_i, \varpi_i}^{-1}, 
\Delta_{\varpi_i, w_0 \varpi_i}^{-1}]_{i \in [1,r]} = 
R_q[G^{w_0,w_0}] = \upClAlg(\sd),
    \]
    where the first equality follows from the definition \eqref{qdBr} of the quantum double Bruhat cells $R_q[G^{w_0, w_0}]$ 
    (note that $I_{w_0,w_0}= 0$). Proposition \ref{prop:Guv_in_bU} implies that $R_q[G] \subset \bUpClAlg$.

    Now the statement of the theorem follows from Theorem C.
\end{proof}

\appendix

\section{Vanishing of non-essential frozen variables}\label{sec:surjection}
In this appendix we address the behavior of partial compactifications with respect to subsets of frozen variables. 

Let $\sd$ denote any classical or quantum seed with the set of vertices $I=I_\ufv\sqcup I_\fv$. Recall that $\nu_i(\ )$ denotes the order of vanishing at $x_i=0$, see Section \ref{sec:cl}. 

Choose and fix a subset $I_\fv'\subset I_\fv$. The partial compactifications associated to $I_\fv'$ are $\bLP^{I_\fv'}:=\bLP^{I_\fv'}(\sd):=\{z\in \LP(\sd)|\nu_i(z)\geq 0,\ \forall i\in I_\fv'\}$ and $\bUpClAlg^{I_\fv'}:=\bUpClAlg^{I_\fv'}(\sd):=\{z\in \upClAlg(\sd)|\nu_i(z)\geq 0,\ \forall i\in I_\fv'\}$. We will omit the superscript ${I_\fv'}$ for simplicity and simply denote them by $\bLP$ and $\bUpClAlg$.

Choose any frozen vertex $j \in I_\fv'$. Note that $\bLP$ is a $\kk[x_j]$-algebra and $\nu_j(\bLP)\geq 0$. Let $A$, $B$ be any $\kk[x_j]$-submodules of $\bLP$. We further require that $A$ is $x_j$-saturated, i.e., $A\cap x_j \bLP=x_j A$. For example, $\bUpClAlg$ is $x_j$-saturated.
\begin{lem}\label{lem:inclusion-quot-alg}
    If $A\subset B$, then this inclusion induces a natural inclusion of the quotient module $A/(x_j A)$ in $B/(x_j B)$.
\end{lem}
\begin{proof}
It suffices to show $A\cap x_j B=x_j A$. Since $A\subset B$, we have $x_j A\subset A\cap x_j B$.  Conversely, since $A$ is $x_j$-saturated, we have $A\cap x_j B\subset A\cap x_j\bLP=x_j A$.
\end{proof}

Choose and fix a subset $F\subset I_\ufv$. Consider the seed $\frz_F\sd$ obtained from $\sd$ by setting the vertices in $F$ to be frozen. Let $\bLP(\frz_F\sd)$ and $\bUpClAlg(\frz_F \sd)$ denote the partial compactification associated to $I_\fv'$. Then we have a natural identification $\bLP=\bLP(\frz_F\sd)$, and a natural inclusion $\bUpClAlg=\cap_{\sd''\in \Delta^+_\sd}\bLP(\sd'') \subset \cap_{\sd''\in \Delta^+_{\frz_F\sd}}\bLP(\sd'') =\bUpClAlg(\frz_F \sd)$

Let $\sd'$ denote the seed obtained from $\frz_F\sd$ by removing the vertex $j\in I_\fv'$. Let $\bLP':=\bLP(\sd')$ and $\bUpClAlg':=\bUpClAlg(\sd')$ denote the partial compactification associated with $I_\fv'\backslash \{j\}$. 

\begin{defn}
    A frozen vertex $j$ is called non-essential (or disconnected) if $b_{jk}=0$ for any unfrozen vertex $k$.
\end{defn}

From now on, we assume that $j$ is non-essential in $\frz_F\sd$. Then we have the canonical surjection $\pi_j:\bLP\rightarrow \bLP'$ sending $x_j$ to $0$.
\begin{prop}\label{prop:quotient_frozen}
    Under the assumption that $j$ is non-essential in $\frz_F\sd$, the map $\pi_j$ restricts to a morphism from $\bUpClAlg(\sd)$ to $\bUpClAlg(\sd')$ whose kernel is the two-sided ideal $x_j \bUpClAlg(\sd)$.
\end{prop}
\begin{proof}

Let $\seq$ denote any sequence of mutations on the set of vertices $I_\ufv\backslash F$. Since $b_{ji}=0$ for all $i\in I_\ufv \backslash F$, $x_j$ does not appear in the exchange relations for computing $x_k(\seq \sd)$, $k\in I_\ufv\backslash F$, from the initial variables in $\sd$. It follows that $\pi_j$ sends $x_k(\seq \sd)\in \bLP$ to $x_k(\seq \sd')\in \bLP'$, inducing an identification $\bLP(\seq \sd)/(x_j \bLP(\seq \sd))\simeq \bLP(\seq \sd')$. Therefore, we obtain $\cap_{\seq}  \bLP(\seq \sd)/(x_j \bLP(\seq \sd)) = \cap_{\seq} \bLP(\seq \sd')=\bUpClAlg'$.

For every such $\seq$, the inclusion $\bUpClAlg\subset \bLP(\seq \sd)$ induces $\bUpClAlg/(x_j \bUpClAlg)\subset \bLP(\seq \sd)/(x_j \bLP(\seq \sd))$ by Lemma \ref{lem:inclusion-quot-alg}. We deduce $\bUpClAlg/(x_j\bUpClAlg) \subset \cap_{\seq}  \bLP(\seq \sd)/(x_j \bLP(\seq \sd)) =\bUpClAlg'$ as desired.
\end{proof}

\section{A geometric proof of the upper cluster algebra structures on classical coordinate rings}

We will take $\kk=\C$ in this section.

We first recall the cluster structure on classical double Bruhat cells following Section \ref{sec:dBruhat}. For any given Weyl group elements $w$ and $u$, choose their reduced words $\uzeta$ and $\ueta$, respectively. Denote $\ubi:=(\uzeta,-\ueta)$ and $\sd^{\BZ}=\dsd(\ubi)^\op$. We have the following isomorphism of commutative algebras, see \eqref{eq:BZ-seed} \eqref{eq:BZ-cluster-var}.
\begin{align}\label{eq:BZ-classical-cluster}
\dkappa^{\BZ}\iota^\op:\upClAlg(\sd^{\BZ}) \simeq \C[G^{u,w}].    
\end{align}
By \cite{shen2021cluster}, different choices of reduced words provide seeds of the same cluster algebra related by mutation sequences. 

Consider the case $w=u=w_0$. Denote $l=l(w_0)$. Let us denote $\sd:=\sd^{\BZ}$, and $\upClAlg=\upClAlg(\sd)$. Identify $\C[G^{w_0,w_0}]$ with $\upClAlg$, such that the initial cluster variables are identified with generalized minors. We will prove the following result.
\begin{thm}\label{thm:classical-cluster-G}
    We have $\bUpClAlg=\C[G].$
\end{thm}

Take any frozen vertex $j$ of $\sd$. Denote $i=|\bi_j|$. Then either $j=-i$ or $j=i^{\max}$. Recall that, in Section \ref{sec:cluster-G}, we have constructed a new seed $\sd_j$ and the associated codimension-$1$ double Bruhat cell $V_j$. Moreover, $j$ is non-essential in $\frz_F \sd_j$, where $F=\{j[1]\}$ in the first case and $F=\{j[-1]\}$ in the second case. And the seed $\sd'_j$ is obtained from $\frz_F\sd_j$ by removing $j$.

In both cases, by Appendix \ref{sec:surjection}, we have a homomorphism associated with the non-essential frozen vertex $j$ in $\frz_F\sd_j$: 
$$ \pi_j: \bUpClAlg(\frz_F\sd_j)\rightarrow \bUpClAlg(\sd'_j),$$
such that $\pi_j(x_j)=0$ and $\pi_j(x_k(\frz_F\sd_j))=x_k(\sd'_j)$ for $k\neq j$. Moreover, we have a cluster structure on $\C[V_j]$ following \eqref{eq:BZ-classical-cluster}:
$$\upClAlg(\sd'_j) \simeq \C[V_j].$$
such that $x_k(\sd'_j)$ corresponds to the restriction of the rational function $x_k(\sd_j)$ on $V_j$. By identifying $x_k(\sd'_j)$ with the restriction of $x_k$ on $V_j$, we will identify $\upClAlg(\sd'_j)$ with $\C[V_j]$ from now on.

\begin{lem}\label{lem:G-incluster-alg}
We have $\C[G]\subset \bUpClAlg$.
\end{lem}
\begin{proof} The frozen variables 
\[
\{\Delta_{w_0 \varpi_i, \varpi_i}, \Delta_{\varpi_i, w_0 \varpi_i} \mid i \in [1,r] \}
\]
are prime elements of $\C[G]$. This follows from the well-known fact that they generate the vanishing ideals of the closures of Bruhat cells $\overline{B^+ w_0 s_i B^+}$ and $\overline{B^- w_0 s_i B^-}$ in $G$, which are irreducible subvarieties since $B^\pm$ are connected algebraic groups. 

    Then the previous arguments in the proofs of Lemma \ref{lem:nu-j-positive} and Proposition \ref{prop:Guv_in_bU} are still valid in the classical case. The desired claim follows.
\end{proof}

\begin{lem}\label{lem:cluster-alg-in-G}
  When $\kk=\C$, we have $\bUpClAlg\subset \C[G]$
\end{lem}
\begin{proof}
Take any $z\in \bUpClAlg$. We identify $\upClAlg$ with $\C[G^{w_0,w_0}]$ as before. Then $z$ is a regular function on $G^{w_0,w_0}$. Recall that $\C[G]$ consists of the functions in $\C[G^{w_0,w_0}]$ which are regular on $G$. By the double Bruhat decomposition, $G\backslash(G^{w_0,w_0}\cup (\cup_{j\in I_\fv} V_j))$ has codimension $2$ in $G$. Since $G$ is smooth and thus normal, a rational function regular at all codimension-$1$ points extends globally. So it suffices to show that the restriction of the rational function $z$ on any $V_j$ is a regular function.

For any given frozen vertex $j$, choose the associated seed $\sd_j$ as before. Note that $z\in \bLP(\sd_j)$ since $z\in \bUpClAlg$. Consider the reduced Laurent expansion $f_z$ of $z$ in $\bLP(\sd_j)=\bLP(\frz_F\sd_j)$, where $\bLP(\frz_F\sd_j):=\{u\in \LP(\frz_F\sd_j)|\nu_i(u)\geq 0, \forall i\in I_{\fv}(\sd_j)\}$. Since $z\in \bUpClAlg$, we have $\nu_j(z)\geq 0$, and $x_j$ does not appear in the denominator of $f_z$. 

We take the Laurent polynomial $f_z|_{x_j\mapsto 0}$ in $\bLP(\sd'_j)=\pi_j\bLP(\frz_F\sd_j)$. Note that it can be viewed as a rational function on $V_j$, and it equals the restriction of the rational function $z$ on $V_j$.

Since $x_j\in \ker \pi_j$, $f_z|_{x_j\mapsto 0}$ is a Laurent expansion of $\pi_j z$ in $\bLP(\sd'_j)=\pi_j \bLP(\frz_F\sd_j)$. Note that $\pi_j z$ is regular on $V_j$, since $\pi_j z\in \bUpClAlg(\sd'_j)\subset \upClAlg(\sd'_j)= \C[V_j]$. Therefore, the rational function $f_z|_{x_j\mapsto 0}$ equals the regular function $\pi_j z$ on $V_j$. So the restriction of the rational function $z$ on $V_j$ is regular.
\end{proof}

\begin{proof}[Proof of Theorem \ref{thm:classical-cluster-G}]
    The desired claim follows from Lemmas \ref{lem:G-incluster-alg} and \ref{lem:cluster-alg-in-G}.
\end{proof}

\def\cprime{$'$}
\providecommand{\bysame}{\leavevmode\hbox to3em{\hrulefill}\thinspace}
\providecommand{\MR}{\relax\ifhmode\unskip\space\fi MR }
\providecommand{\MRhref}[2]{%
  \href{http://www.ams.org/mathscinet-getitem?mr=#1}{#2}
}
\providecommand{\href}[2]{#2}


\end{document}